\documentclass{article}
\usepackage{amssymb,named}

\newcommand{\Bcite}[2]{\cite[#1]{#2}}
\newcommand{\pp}[1]{#1}

\let\epsilon\varepsilon

\newbox\gnBoxA
\newdimen\gnCornerHgt
\setbox\gnBoxA=\hbox{$\ulcorner$}
\global\gnCornerHgt=\ht\gnBoxA
\newdimen\gnArgHgt

\def\Godelnum #1{%
        \setbox\gnBoxA=\hbox{$#1$}%
        \gnArgHgt=\ht\gnBoxA%
        \ifnum \gnArgHgt<\gnCornerHgt
                \gnArgHgt=0pt%
        \else
                \advance \gnArgHgt by -\gnCornerHgt%
        \fi
        \raise\gnArgHgt\hbox{$\ulcorner$} \box\gnBoxA %
                \raise\gnArgHgt\hbox{$\urcorner$}}

\title{Hilbert's Program Then and Now}

\author{Richard Zach}
\date{}
\begin{document}\bibliographystyle{named}

\maketitle

\begin{abstract}
Hilbert's program was an ambitious and wide-ranging project in the philosophy
and foundations of mathematics.  In order to ``dispose of the foundational
questions in mathematics once and for all,'' Hilbert proposed a two-pronged
approach in 1921: first, classical mathematics should be formalized in
axiomatic systems; second, using only restricted, ``finitary'' means, one
should give proofs of the consistency of these axiomatic systems.  Although
G\"odel's incompleteness theorems show that the program as originally
conceived cannot be carried out, it had many partial successes, and generated
important advances in logical theory and meta-theory, both at the time and
since.  The article discusses the historical background and development of
Hilbert's program, its philosophical underpinnings and consequences, and its
subsequent development and influences since the 1930s.

\end{abstract}

\section{Introduction}

Hilbert's program\index{Hilbert's program|(} is, in the first instance, a
proposal and a research program in the philosophy and foundations of
mathematics.  It was formulated in the early 1920s by German mathematician
David Hilbert\index{Hilbert, D.|(} (1862--1943), and was pursued by him and his
collaborators at the University of G\"ottingen and elsewhere in the 1920s and
1930s.  Briefly, Hilbert's proposal called for a new foundation of mathematics
based on two pillars: the axiomatic method\index{axiomatic method}, and
finitary proof theory.  Hilbert thought that by formalizing mathematics in
axiomatic systems, and subsequently proving by finitary methods that these
systems are consistent (i.e., do not prove contradictions), he could provide a
philosophically satisfactory grounding of classical, infinitary mathematics
(analysis and set theory).  Had it been successful, Hilbert's program would
perhaps not only have eclipsed in subsequent influence other foundational
projects of the time, such as the logicist projects pursued by Frege and
Russell (see Vol.~9, Ch.~5) and Brouwer's\index{Brouwer, L. E. J.}
intuitionism\index{intuitionism} (see \textbf{Handbook reference?}), but it
would also have achieved Hilbert's stated goal, viz., to ``dispose of the
foundational questions in mathematics as such once and for all''
\cite[228]{Hilbert:29}.  Unfortunately, G\"odel's theorems (see Ch.~11) show
that the program as originally envisaged by Hilbert cannot be carried out.

Although Hilbert's own project for the foundation of mathematics was
ultimately unsuccessful, the project itself and technical advances made in its
pursuit have had an enormous impact on logic and the foundations of
mathematics more generally.  In order to carry out the first part of the
program, the axiomatization of mathematics in a formal system, Hilbert and his
collaborators had pushed forward the development of logical formalisms in
which such an axiomatization could be carried out.  This led, in particular,
to the first axiomatizations of propositional and first-order logic as
independent systems (i.e., other than as fragments of more comprehensive
higher-order systems, such as Frege's\index{Frege, G.} \emph{Begriffsschrift}
and the Whitehead-Russell\index{Russell, B.}\index{Whitehead, A. N.} system of
\emph{Principia mathematica})\index{Principia Mathematica} and their
meta-logical investigation, such as the proof (by Paul Bernays\index{Bernays,
  P.} in 1918) of the completeness of the propositional calculus and the work
in Hilbert's school from 1921 onward on the decision problem.  The
investigation of the meta-logical properties of logical systems led directly
to some of the most important meta-logical results in logic, viz., G\"odel's
completeness theorem\index{completeness theorem} and the negative solution by
Church\index{Church, A.} and Turing\index{Turing, A. M.} of the decision
problem\index{decision problem}.  The development of proof theory itself is an
outgrowth of Hilbert's program.  Gentzen's\index{Gentzen, G.} development of
natural deduction and the sequent calculus was carried out in the tradition of
Hilbert's program and with the aim of constructing a logical system which
facilitates consistency proofs\index{consistency proof}.
G\"odel\index{G\"odel, K.} obtained his incompleteness
theorems\index{incompleteness theorems} while trying to prove the consistency
of analysis.  And the tradition of reductive proof theory of the
Gentzen-Sch\"utte school and others is itself a direct continuation of
Hilbert's program.

The present chapter is divided into three parts: The first part provides a
sketch of the historical development of logic, proof theory, and philosophy of
mathematics in the work of Hilbert and his followers through the 1930s.  The
second part deals with the philosophical interpretation and assessment of
Hilbert's program.  The third part presents recent work in proof theory which
bears on the aims of Hilbert's program.

\section{Hilbert's program then}

\subsection{Hilbert's early work on foundations}

Hilbert's\index{Hilbert's program!history of|(} work on the foundations of
mathematics can be traced to his work on geometry of the 1890s which resulted
in his influential textbook \emph{Foundations of Geometry}
\shortcite{Hilbert:99}.  One philosophical advance of this work was the
development of Hilbert's conception of the axiomatic method\index{axiomatic
  method}.  Hilbert believed that the proper way to develop any scientific
subject rigorously required an axiomatic approach.  In providing an axiomatic
treatment, the theory would be developed independently of any need for
intuition, and it would facilitate an analysis of the logical relationships
between the basic concepts and the axioms.  Of basic importance for an
axiomatic treatment are, so Hilbert, investigation of the independence and,
above all, of the consistency of the axioms. In his 1902 lectures on the
foundations of geometry, he puts it thus:
\begin{quotation}
  Every science takes its starting point from a sufficiently coherent body of
  facts is given. It takes form, however, only by \emph{organizing} this body
  of facts. This organization takes place through the \emph{axiomatic
    method}\index{axiomatic method}, i.e., one constructs a \emph{logical
    structure of concepts} so that the relationships between the concepts
  correspond to relationships between the facts to be organized.

There is arbitrariness in the construction of such a structure of concepts;
we, however, \emph{demand} of it:

1) completeness, 2) independence, 3) consistency. \cite[540]{Hilbert:04}
\end{quotation}
From the time of his work on geometry forward, the last consideration,
consistency, was of special importance in Hilbert's conception of the
axiomatic method\index{axiomatic method} in general and the foundations of
mathematics in particular.  Hilbert was heavily influenced by the foundational
views of late-19th century mathematicians, in particular, Cantor\index{Cantor,
  G.}, Dedekind,\index{Dedekind, R.} and Kronecker\index{Kronecker, L.|(}.  He
shared with Dedekind and Cantor the view that mathematical activity should be
free of constraints, which led to his view, highlighted in his correspondence
with Frege\index{Frege, G.}, that consistency of an axiomatic theory
guarantees the existence of the structure described, and is in this sense
sufficient to justify the use of the theory.  And he shared with Kronecker a
recognition that elementary arithmetic has a privileged role in mathematics,
although he was of course opposed to the converse espoused by Kronecker, viz.,
that the natural numbers, and constructions based on elementary arithmetic,
exhaust legitimate mathematics.  These two influences in Hilbert's thought are
at the root of his investigations of consistency.

Proofs of consistency\index{consistency proof} for the axioms of geometry can
be given by providing an interpretation of the system in the real plane, and
thus the consistency of geometry is reduced to the consistency of analysis.
Analysis, of course, itself requires justification.  In
\shortcite{Hilbert:00}, Hilbert approached the problem from the axiomatic
standpoint by proposing an axiomatization of the real numbers.  In order to
show the consistency of this system, Hilbert expressly rejected the
construction of a model, e.g., a construction based on Dedekind
cuts\index{Dedekind cut} of rationals, as an option. He considered the
construction of the reals from the rationals and ultimately the natural
numbers using the ``genetic method'' as insufficient: ``Despite the high
pedagogical and heuristic value of the genetic method, for the final
presentation and the complete logical grounding of our knowledge the axiomatic
method deserves the first rank'' \cite[1093]{Hilbert:00}.  Hilbert thus was
after a \emph{direct} consistency proof\index{consistency proof} of analysis,
i.e., one not based on reduction to another theory.  He proposed the problem
of finding such a proof as the second of his 23 mathematical problems in his
address to the International Congress of Mathematicians in 1900
\shortcite{Hilbert:00a}.

The discovery of Russell's paradox\index{Russell's paradox} in 1902 made it
clear that an axiomatic foundation of arithmetic and set theory requires a
more precise development of the underlying logical systems.  Hilbert knew of
the paradoxes of set theory\index{set theory!paradoxes of} from
Cantor\index{Cantor, G.} and Zermelo\index{Zermelo, E.|(}, but it was
apparently not until Russell\index{Russell, B.}'s \shortcite{Russell:02a}
publication of the contradiction in Frege\index{Frege, G.}'s system that
Hilbert and Zermelo realized their importance.  Hilbert's exchange with
Frege\index{Frege, G.} on the axiomatic approach to geometry led him to
realize that his conceptions of ``axiom,'' ``definition,'' ``proof'' were in
need of clarification.  In response, Hilbert formulated an early version of
his proof-theoretical program in his 1904 Heidelberg talk
\shortcite{Hilbert:05}.  After criticizing the foundational views of
Kronecker\index{Kronecker, L.|)} as dogmatic, and those of Frege\index{Frege,
  G.} and Dedekind\index{Dedekind, R.} as suffering from ``unavoidable
contradictions,'' he writes:
\begin{quotation}
  Arithmetic is often considered to be a part of logic, and the traditional
  fundamental logical notions are usually presupposed when it is a question of
  establishing a foundation for arithmetic. If we observe attentively,
  however, we realize that in the traditional exposition of the laws of logic
  certain fundamental arithmetic notions are already used, for example, the
  notion of set and, to some extent, also that of number. Thus we find
  ourselves turning in a circle, and that is why a partly simultaneous
  development of the laws of logic and of arithmetic is required if paradoxes
  are to be avoided. \cite[131]{Hilbert:05}
\end{quotation}
Hilbert's sketch of this ``simultaneous development'' of logic and arithmetic
in the case of a very basic axiom system for the natural numbers is very close
to the the approach Hilbert's proof theoretic program would take 20 years
later: Hilbert gives a direct argument that no contradiction can arise from
the five axioms of his system.  

This was a promising start, but several factors delayed the further
development of Hilbert's proof theoretic program.  One was Poincar\'e's
\shortcite{Poincare:06a}\index{Poincar\'e, H.} criticism of what he saw as a
viciously circular use of induction\index{induction!Poincar\'e on} in
Hilbert's sketched consistency proof\index{consistency proof} (see
\cite{Steiner:75}\index{Steiner, M.}, Appendix).  Moreover, Hilbert realized
that axiomatic investigations required a well worked-out logical formalism in
which axiomatic systems could be developed.  At the time he used a logical
formalism based on Schr\"oder's\index{Schr\"oder, E.}  algebraic
logic\index{logic!algebraic}, which was not particularly suited as a formalism
for the axiomatization of mathematics.  Following the 1905 lectures on
foundations, Hilbert turned his immediate attention to work in other areas of
mathematics and theoretical physics.  He did, however, actively support others
who worked on foundational questions in G\"ottingen, in particular Ernst
Zermelo\index{Zermelo, E.|)} and Leonard Nelson\index{Nelson,
  L.}.\footnote{See \cite{Sieg:99,Sieg:02}\index{Sieg, W.},
  \cite{Stein:88}\index{Stein, H.},
  \cite{Hallett:90,Hallett:94}\index{Hallett, M.},
  \cite{Mancosu:98c}\index{Mancosu, P.}, and \cite{Avigad:01}\index{Avigad,
    J.}\index{Reck, E.} for further discussion of the conceptual framework and
  historical background of Hilbert's thought, and
  \cite{Resnik:74}\index{Resnik, M.} on the Frege\index{Frege, G.}-Hilbert
  correspondence. On Hilbert's foundational interests before 1917, and his
  engagement for Husserl\index{Husserl, E.}, Zermelo, and Nelson\index{Nelson,
    L.} in G\"ottingen, see \cite{Peckhaus:90}\index{Peckhaus, V.}.  On
  general discussions of formalism and the place of Hilbert's thought in the
  mathematical context of the late 19th century, see
  \cite{Webb:97}\index{Webb, J.} and \cite{Detlefsen:05}\index{Detlefsen,
    M.}.}

The publication of Whitehead\index{Whitehead, A. N.} and Russell\index{Russell, B.}'s
\emph{Principia Mathematica}\index{Principia Mathematica}
\shortcite{Russell:10,Russell:12,Russell:13} provided the required logical
basis for a renewed attack on foundational issues.  Beginning in 1914,
Hilbert's student Heinrich Behmann\index{Behmann, H.} and others studied the
system of \emph{Principia}.\footnote{See \cite{Mancosu:99} and
  \shortcite{Mancosu:03}\index{Mancosu, P.} on Behmann\index{Behmann, H.}'s
  role in Hilbert's school and the influence of Russell\index{Russell, B.}.} Hilbert
himself returned to work on the foundations of mathematics in 1917.  In
September 1917, he delivered an address to the Swiss Mathematical Society
entitled ``Axiomatic Thought'' \shortcite{Hilbert:18}.  It is his first
published contribution to mathematical foundations since 1905.  In it, he
again emphasized the requirement of consistency proofs\index{consistency
  proof} for axiomatic systems: ``The chief requirement of the theory of
axioms must go farther [than merely avoiding known paradoxes], namely, to show
that within every field of knowledge contradictions based on the underlying
axiom-system are \emph{absolutely impossible}.'' He posed the proof of the
consistency of the integers (and of set theory) again as the main problems.
In both these cases, there seems to be nothing more fundamental available to
which the consistency could be reduced other than logic itself.  Hilbert at
the time considered the problem as essentially solved by
Whitehead\index{Whitehead, A. N.} and Russell\index{Russell, B.}'s work in
\emph{Principia}\index{Principia Mathematica}.  Nevertheless, other
fundamental problems of axiomatics remained unsolved, including the problem of
the ``decidability of every mathematical question,'' which also traces back to
Hilbert's 1900 address.

These unresolved problems of axiomatics led Hilbert to devote significant
effort to work on logic in the following years.  In 1917, Paul
Bernays\index{Bernays, P.} joined him as his assistant in G\"ottingen.  In a
series of courses from 1917--1921, Hilbert, with the assistance of
Bernays\index{Bernays, P.} and Behmann\index{Behmann, H.}, made significant
new contributions to formal logic.  The course from 1917 \cite{Hilbert:17}, in
particular, contains a sophisticated development of first-order
logic\index{logic!first order}, and forms the basis of Hilbert and
Ackermann\index{Ackermann, W.}'s textbook \emph{Principles of Theoretical
  Logic} \shortcite{Hilbert:28b}.  In 1918, Bernays\index{Bernays, P.}
submitted a treatise on the propositional calculus of \emph{Principia
  mathematica} as a \emph{Habilitationsschrift}; it contains the first
completeness proof\index{completeness!of propositional logic} of the
propositional calculus for truth-functional semantics.\footnote{See
  \cite{Moore:97}\index{Moore, G.}, \cite{Sieg:99}\index{Sieg, W.} and
  \cite{Zach:99}\index{Zach, R.} for more detail on the
  development of logic in Hilbert's school around 1918.}

The 1917--18 lectures were only the beginning of a strand of work on logic and
meta-logic in Hilbert's school, including work on the decision problem.  The
decision problem for first-order logic was tightly bound up with the aim of
finding a completeness proof for the first-order predicate
calculus\index{logic!first order} (the ``restricted calculus of functions'' in
Hilbert's terminology).  This aim was stated in the 1917--18 lectures, but
since completeness does not hold for first-order logic in any purely syntactic
sense (an early result due to Ackermann\index{Ackermann, W.}), a development
of the semantics of first-order logic was needed first. The decision
problem\index{decision problem}, one of Hilbert's main aims for
metamathematics in the 1920s, was already an issue in the lectures from 1905,
and has its roots in Hilbert's belief, first explicitly stated in the Paris
address, that ``in mathematics, there is no ignorabimus,'' i.e., that every
mathematical question can be solved either affirmatively or negatively.  The
questions of completeness\index{completeness} and decidability thus became
closely linked in the 1920s, with Behmann\index{Behmann, H.},
Bernays\index{Bernays, P.}, Sch\"onfinkel, and later
Ackermann\index{Ackermann, W.} working on special cases of the decision
problem\index{decision problem} for first-order logic throughout the 1920s.

\subsection{The consistency program, finitism, and proof theory}

In about 1920, Hilbert came to reject Russell\index{Russell, B.}'s
logicist\index{logicism} solution to the consistency problem for arithmetic,
mainly for the reason that the axiom of reducibility cannot be accepted as a
purely logical axiom.  In lectures from the Summer term 1920, he concluded
that ``the aim of reducing set theory, and with it the usual methods of
analysis, to logic, has not been achieved today and maybe cannot be achieved
at all'' \cite{Hilbert:20b}. At the same time, Brouwer's\index{Brouwer, L. E.
  J.} intuitionist\index{intuitionism} mathematics gained currency.  In
particular, Hilbert's former student Hermann Weyl\index{Weyl, H.}  converted
to intuitionism.  Weyl's 1920 address to the Hamburg Mathematical Seminar,
``The new foundational crisis in mathematics'' \shortcite{Weyl:21} was
answered by Hilbert in three talks in Hamburg in the Summer of 1921
\shortcite{Hilbert:22a}.  Here, Hilbert presented his own mature proposal for
a solution to the problem of the foundation of mathematics.  This proposal
incorporated Hilbert's ideas from 1904 regarding direct consistency
proofs\index{consistency proof}, his conception of axiomatic systems, and also
the technical developments in the axiomatization of mathematics in the work of
Russell\index{Russell, B.} as well as the further developments carried out by him and
his collaborators.  What was new was the way in which Hilbert wanted to imbue
his consistency project with the philosophical significance necessary to
answer Brouwer\index{Brouwer, L. E. J.} and Weyl\index{Weyl, H.}'s criticisms:
the finitary point of view.

According to Hilbert, there is a privileged part of
mathematics,\index{finitism|(} contentual elementary number theory, which
relies only on a ``purely intuitive basis of concrete signs.''  Whereas the
operating with abstract concepts was considered ``inadequate and uncertain,''
there is a realm of
\begin{quote}
  extra-logical discrete objects, which exist intuitively as immediate
  experience before all thought.  If logical inference is to be certain, then
  these objects must be capable of being completely surveyed in all their
  parts, and their presentation, their difference, their succession (like the
  objects themselves) must exist for us immediately, intuitively, as something
  which cannot be reduced to something else.\footnote{\cite[202]{Hilbert:22a}.
    The passage is repeated almost verbatim in \cite[376]{Hilbert:26},
    \cite[464]{Hilbert:28}, and \cite[267]{Hilbert:31a}}
\end{quote}
The objects in questions are \emph{signs}, both numerals and the signs
that make up formulas a formal proofs.  The domain of contentual number theory
consists in the finitary numerals, i.e., sequences of strokes.  These have no
meaning, i.e., they do not stand for abstract objects, but they can be
operated on (e.g., concatenated) and compared.  Knowledge of their properties
and relations is intuitive and unmediated by logical inference. Contentual
number theory developed this way is secure, according to Hilbert: no
contradictions can arise simply because there is no logical structure in the
propositions of contentual number theory.\index{finitism|)}

The intuitive-contentual operations with signs form the basis of Hilbert's
meta-mathematics\index{meta-mathematics}.  Just as contentual number theory
operates with sequences of strokes, so meta-mathematics operates with
sequences of symbols (formulas, proofs).  Formulas and proofs can be
syntactically manipulated, and the properties and relationships of formulas
and proofs are similarly based in a logic-free intuitive capacity which
guarantees certainty of knowledge about formulas and proofs arrived at by such
syntactic operations.  Mathematics itself, however, operates with abstract
concepts, e.g., quantifiers, sets, functions, and uses logical inference based
on principles such as mathematical induction or the principle of the excluded
middle.  These ``concept-formations'' and modes of reasoning had been
criticized by Brouwer and others on grounds that they presuppose infinite
totalities as given, or that they involve impredicative\index{impredicative}
definitions (which were considered by the critics as viciously circular).
Hilbert's aim was to justify their use.  To this end, he pointed out that they
can be formalized in axiomatic systems (such as that of
\emph{Principia}\index{Principia Mathematica} or those developed by Hilbert
himself), and mathematical propositions and proofs thus turn into formulas and
derivations from axioms according to strictly circumscribed rules of
derivation.  Mathematics, so Hilbert, ``becomes an inventory of provable
formulas.'' In this way the proofs of mathematics are subject to
metamathematical, contentual investigation.  The goal of Hilbert's program is
then to give a contentual, meta-mathematical proof that there can be no
derivation of a contradiction\index{contradiction}, i.e., no formal derivation
of a formula $A$ and of its negation $\neg A$.

This sketch of the aims of the program was fleshed out by Hilbert and his
collaborators in the following 10 years.  On the conceptual side, the finite
standpoint and the strategy for a consistency proof\index{consistency proof}
were elaborated by Hilbert \shortcite{Hilbert:23,Hilbert:26,Hilbert:28} and
Bernays\index{Bernays, P.}  \shortcite{Bernays:22,Bernays:28a,Bernays:30}, of
which Hilbert's article ``On the infinite'' \shortcite{Hilbert:26} provides
the most detailed discussion of the finitary standpoint.  In addition to
Hilbert and Bernays\index{Bernays, P.}, a number of other people were involved
in technical work on the program.  The $\epsilon$-operator\index{epsilon
  calculus@$\epsilon$-calculus|(} was first introduced in the Hamburg lecture
of 1921 \cite{Hilbert:22a}, and developed in a number of lectures given in
G\"ottingen \cite{Hilbert:21,Hilbert:22b}, as well as in \cite{Hilbert:23}.
Hilbert and Bernays\index{Bernays, P.}  developed the $\epsilon$-calculus as
their definitive formalism for axiom systems for arithmetic and analysis, and
the so-called $\epsilon$-substitution method as the preferred approach to
giving consistency proofs.

Briefly, the $\epsilon$-calculus is a formalism that includes $\epsilon$ as a
term-forming operator.  If $A(x)$ is a formula, then $\epsilon_x A(x)$ is a
term, which intuitively stands for a witness for~$A(x)$.  In a logical
formalism containing the $\epsilon$-operator, the quantifiers can be defined
by: $\exists x\, A(x) \equiv A(\epsilon_x A(x))$ and $\forall x\, A(x) \equiv
A(\epsilon_x \lnot A(x))$.  The only additional axiom necessary is the
so-called ``transfinite axiom,'' $A(t) \to A(\epsilon_x A(x))$.  Based on this
idea, Hilbert and his collaborators developed axiomatizations of number theory
and analysis.  Consistency proofs\index{consistency proof} for these systems
were then given using the $\epsilon$-substitution method\index{epsilon
  substitution@$\epsilon$-substitution method}.  The idea of this method is,
roughly, that the $\epsilon$-terms $\epsilon_x A(x)$ occurring in a formal
proof are replaced by actual numerals, resulting in a quantifier-free proof.
The simplest case, outlined in Hilbert's papers, is as follows.  Suppose we
had a (suitably normalized) derivation of $0 = 1$ that contains only one
$\epsilon$-term $\epsilon_x A(x)$.  Replace all occurrences of $\epsilon_x
A(x)$ by $0$.  The instances of the transfinite axiom then are all of the form
$A(t) \to A(0)$.  Since no other $\epsilon$-terms occur in the proof, $A(t)$
and $A(0)$ are basic numerical formulas without quantifiers and, we may
assume, also without free variables.  So they can be evaluated by finitary
calculation.  If all such instances turn out to be true numerical formulas, we
are done.  If not, this must be because $A(t)$ is true for some~$t$, and
$A(0)$ is false.  Then replace $\epsilon_x A(x)$ instead by $n$, where $n$ is
the numerical value of the term~$t$.  The resulting proof is then seen to be a
derivation of $0=1$ from true, purely numerical formulas using only modus
ponens, and this is impossible.  Indeed, the procedure works with only slight
modifications even in the presence of the induction axiom, which in the
$\epsilon$-calculus takes the form of a least number principle: $A(t) \to
\epsilon_x A(x) \le t$, which intuitively requires $\epsilon_x A(x)$ to be the
\emph{least} witness for $A(x)$.

The $\epsilon$-substitution method is simple enough for the basic cases
considered by Hilbert, but becomes extremely complex when $\epsilon$-operators
are nested.  In his 1924 dissertation, \cite{Ackermann:24}\index{Ackermann,
  W.} presented an (erroneous) consistency proof\index{consistency proof}
based on the $\epsilon$-substitution method for a version of analysis.  John
von Neumann\index{von Neumann, J.}, then visiting G\"ottingen, gave a
corrected consistency proof\index{consistency proof} for a system of the
$\epsilon$-formalism (which, however, did not include the induction axiom) in
1925 \shortcite{Neumann:27}\index{von Neumann, J.}.  Building on von Neumann's
work, Ackermann\index{Ackermann, W.} devised a new $\epsilon$-substitution
procedure which he communicated to Bernays\index{Bernays, P.} (see
\cite{Bernays:28a}).  Ackermann and Bernays considered the proof to be correct
for the entire first-order fragment of arithmetic and were confident that it
could be extended to a consistency proof\index{consistency proof} of analysis.
In his address ``Problems of the grounding of mathematics'' to the
International Congress of Mathematicians in Bologna in 1928
\shortcite{Hilbert:29}, Hilbert optimistically claimed that the work of
Ackermann\index{Ackermann, W.} and von Neumann\index{von Neumann, J.} had
established the consistency of number theory and that the proof for analysis
had already been carried out by Ackermann\index{Ackermann, W.} ``to the extent
that the only remaining task consists in the proof of an elementary finiteness
theorem that is purely arithmetical.''\footnote{See
  \cite{Avigad:02}\index{Avigad, J.} for a general introduction to the
  $\epsilon$-calculus and \cite{Zach:03}\index{Zach, R} and
  \shortcite{Zach:02} on the history of the $\epsilon$-calculus and the
  substitution method. \cite{Sieg:99}\index{Sieg, W.} presents a detailed and
  perceptive analysis of the development of the program and its influence as a
  whole.}\index{epsilon calculus@$\epsilon$-calculus|)}

\subsection{The impact of G\"odel's incompleteness 
theorems}\label{impact-godel}

G\"odel\index{G\"odel, K.}'s incompleteness theorems\index{incompleteness
  theorems|(} \cite{Godel:31} showed that Hilbert's optimism was unfounded. In
September 1930, Kurt G\"odel\index{G\"odel, K.} announced his first
incompleteness theorem\index{incompleteness theorems!first incompleteness
  theorem} at a conference in K\"onigsberg.  Von Neumann\index{von Neumann,
  J.}, who was in the audience, immediately recognized the significance of
G\"odel\index{G\"odel, K.}'s result for Hilbert's program.  Shortly after the
conference he wrote to G\"odel\index{G\"odel, K.}, telling him that he had
found a corollary to G\"odel's result. G\"odel\index{G\"odel, K.} had found
the same result already independently: the second incompleteness
theorem\index{incompleteness theorems!second incompleteness theorem},
asserting that the system of \emph{Principia}\index{Principia Mathematica}
does not prove the formalization of the claim that the system of
\emph{Principia} is consistent (provided it is).  All the methods of finitary
reasoning used in the consistency proofs\index{consistency proof} up till then
were believed to be formalizable in \emph{Principia}, however.  Hence, if the
consistency of \emph{Principia} were provable by the methods used in
Ackermann\index{Ackermann, W.}'s proofs, it should be possible to formalize
this proof in \emph{Principia}; but this is what the second incompleteness
theorem\index{incompleteness theorems!second incompleteness theorem} states is
impossible.  Bernays\index{Bernays, P.}  also immediately realized the
importance of G\"odel\index{G\"odel, K.}'s results after he studied G\"odel's
paper in January 1931. He wrote to G\"odel\index{G\"odel, K.} that (under the
assumption that finitary reasoning can be formalized in \emph{Principia}) the
incompleteness theorems seem to show that a finitary consistency
proof\index{consistency proof} of \emph{Principia} is impossible.  Shortly
thereafter, von Neumann\index{von Neumann, J.} showed that
Ackermann\index{Ackermann, W.}'s consistency proof\index{consistency proof} is
flawed and provided a counterexample to the proposed $\epsilon$-substitution
procedure.\footnote{The correspondence between G\"odel\index{G\"odel, K.}  and
  Bernays\index{Bernays, P.} is published in \cite[41--313]{Godel:03} and that
  with von Neumann\index{von Neumann, J.} in \cite[327--377]{Godel:03a}.  See
  also the informative introductions by Feferman and Sieg, respectively, to
  these sections of the correspondence, as well as
  \cite{Mancosu:04a}\index{Mancosu, P.} and the last section of
  \cite{Zach:03}\index{Zach, R.}.}

Although the impact of G\"odel\index{G\"odel, K.}'s incompleteness theorems
for Hilbert's program was recognized soon after their publication, Hilbert's
program was by no means abandoned.  Hilbert himself was no longer actively
involved in foundational research, but Bernays\index{Bernays, P.} continued to
work on the foundations of mathematics.  The two-volume \emph{Grundlagen der
  Mathematik} \cite{HilbertBernays:34,HilbertBernays:39} was prepared by
Bernays\index{Bernays, P.} alone, and included new results by
Ackermann\index{Ackermann, W.} and Bernays\index{Bernays, P.} on the
$\epsilon$-calculus.  It also included Herbrand's\index{Herbrand, J.}
\shortcite{Herbrand:30} work on proof theory, and a sketch of
Gentzen\index{Gentzen, G.}'s \shortcite{Gentzen:36} consistency
proof\index{consistency proof} of first-order arithmetic\index{Peano
  arithmetic}.  Bernays\index{Bernays, P.}'s and Gentzen\index{Gentzen, G.}'s
work, in particular, focused on possible extensions of the finitary standpoint
which could yield consistency proofs for substantial parts of mathematics in
spite of G\"odel\index{G\"odel, K.}'s theorems.\index{incompleteness
  theorems|)}

Gentzen\index{Gentzen, G.}'s first consistency proof\index{consistency proof}
for number theory, the so-called galley proof \shortcite{Gentzen:35}, was the
result of a combination of Gentzen's \shortcite{Gentzen:34} earlier work on
the logical formalism of the sequent calculus, which provided a
proof-theoretically more convenient axiomatization of arithmetic, and a new
strategy of proving consistency.  This strategy involved defining certain
``reduction steps'' on proofs: local transformations of parts of a derivation
in the new formalism.  The consistency proof\index{consistency proof} then
proceeds by showing that by iterating these reductions on a proof one
eventually arrives at a proof of a special form (a proof free of the cut rule
and the induction rule), and no proof of such a form can be the proof of a
contradiction.  The first version of the proof relied on the notion of a
reduction rule, which itself cannot be formalized in arithmetic.\footnote{On
  the galley proof, see \cite{Bernays:70}\index{Bernays, P.}, \cite[Appendix
  II]{Kreisel:71}, and \cite{Negri:80}.}  This notion met with some
objections, and in the revised, published version \shortcite{Gentzen:36},
Gentzen\index{Gentzen, G.} replaced the appeal to reduction rules by a proof
that the iteration of reduction steps itself terminates.  He did this by
assigning a measure for the complexity of a given derivation, and showed that
the result of the application of a reduction step to a proof reduces the
complexity measure of that proof.  The complexity measure
Gentzen\index{Gentzen, G.} used was certain finite strings of numbers which
may be interpreted as naming countable ordinals less
than~$\epsilon_0$\index{epsilon zero@$\epsilon_0$}.\footnote{If $\omega_0 =
  \omega$, and $\omega_{n+1} = \omega^{\omega_n}$, then the ordinal
  $\epsilon_0$ is the limit of $\omega_n$ for $n = 1, 2, \ldots$. In other
  words, $\epsilon_0$ is the least fixed point of $\alpha = \omega^\alpha$.}
The consistency result\index{consistency proof} then follows if one accepts
that there are no infinite descending sequences of such ordinal notations, or,
more precisely, by using transfinite induction up to~$\epsilon_0$.  This
principle, by G\"odel\index{G\"odel, K.}'s incompleteness
theorem\index{incompleteness theorems}, cannot itself be formalized in
first-order arithmetic \cite{Gentzen:43}.  Gentzen\index{Gentzen, G.}'s proof
allowed \cite{Ackermann:40}\index{Ackermann, W.} to give a correct consistency
proof\index{consistency proof} based on the $\epsilon$-substitution method for
first-order arithmetic, also based on transfinite induction up
to~$\epsilon_0$.

Gentzen\index{Gentzen, G.}'s work marks the beginning of post-G\"odelian proof
theory.  In the proof-theoretic tradition of Gentzen\index{Gentzen, G.},
axiomatic theories are analyzed according to which transfinite induction
principles are required to prove the consistency\index{consistency proof} of
the theory.  However, Gentzen\index{Gentzen, G.}'s contribution and influence
goes beyond this: He emphasized that proof-theoretic methods do not just allow
us to prove the consistency of a theory, but that they also enable us to
extract information from proofs beyond the fact that the formula proved
follows from the axioms.\index{Hilbert's program!history of|)}

\section{Philosophical interpretation of Hilbert's program}

The philosophical importance and influence of Hilbert's work on foundations is
twofold.  First, the epistemological standpoint of Hilbert's finitism is of
interest in the philosophy of mathematics quite independently of the success
of the proof-theoretic program which it underlies.  Like other important
proposals in the philosophy of mathematics such as
intuitionism\index{intuitionism}, predicativism, and logicism\index{logicism},
Hilbert's finitism is, \emph{inter alia}, a philosophical view about the
nature of mathematical knowledge and delineates a particular set of
(finitarily) meaningful propositions, as well as finitarily admissible
constructions and methods of proof.  Debate about the status of finitary
evidence and proof are still very much alive today.  Second, Hilbert's program
can and has been seen as a version of reductive instrumentalism in
mathematics.  That is to say, one can read Hilbert as proposing that only a
certain part of mathematics (propositions, proofs) is meaningful, viz., the
finitary part. The rest, which includes classical infinitary mathematics (full
first-order arithmetic\index{Peano arithmetic}, analysis, and set theory, in
particular) are mere instruments.  Hilbert's program has thus been an
important inspiration for related instrumentalist proposals in the philosophy
of mathematics (e.g., \cite{Field:80,Detlefsen:86}\index{Field,
  H.}\index{Detlefsen, M.}).

This section will expand on these two themes.  In the case of both the debate
on the philosophy of finitism, and on the view of Hilbert's program as an
instrumentalist philosophy of mathematics, questions of historical
interpretation interact with conceptual analysis.  The distinction between
these aspects should be kept in mind.

\subsection{The finitary point of view}

The cornerstone of Hilbert's philosophy of mathematics, and the substantially
new aspect of his foundational thought from \shortcite{Hilbert:22a} onward,
was the so-called finitary standpoint\index{finitism|(}.  This methodological
standpoint consists in a restriction of mathematical thought to objects which
are ``intuitively present as immediate experience prior to all thought,'' and
to those operations on and methods of reasoning about such objects which do
not require the introduction of abstract concepts, in particular, require no
appeal to completed infinite totalities.

Hilbert characterized the domain of finitary reasoning in a well-known
paragraph which appears in roughly the same formulation in all of
Hilbert's more philosophical papers from the 1920s
\shortcite{Hilbert:22a,Hilbert:26,Hilbert:28,Hilbert:31a}:
\begin{quote}
[A]s a condition for the use of logical inferences and the performance
of logical operations, something must already be given to our faculty
of representation, certain extra-logical concrete objects that are
intuitively present as immediate experience prior to all thought.  If
logical inference is to be reliable, it must be possible to survey
these objects completely in all their parts, and the fact that they
occur, that they differ from one another, and that they follow each
other, or are concatenated, is immediately given intuitively, together
with the objects, as something that can neither be reduced to anything
else nor requires reduction.  This is the basic philosophical position
that I consider requisite for mathematics and, in general, for all
scientific thinking, understanding, and
communication. \cite[376]{Hilbert:26}
\end{quote}
These objects are, for Hilbert, the \emph{signs}. For the domain of contentual
number theory, the signs in question are sequences of strokes (``numerals'')
such as
\[
\mid, \mid\mid, \mid\mid\mid, \mid\mid\mid\mid\mid.
\]
The question of how exactly Hilbert understood the numerals is difficult to
answer.  What is clear in any case is that they are logically primitive, i.e.,
they are neither concepts (as Frege\index{Frege, G.}'s numbers are) nor sets.
For Hilbert, the important issue is not primarily their metaphysical status
(abstract versus concrete in the current sense of these terms), but that they
do not enter into logical relations, e.g., they cannot be predicated of
anything.  In Bernays\index{Bernays, P.}'s most mature presentations of
finitism \cite{Bernays:30,HilbertBernays:39}, the objects of finitism are
characterized as \emph{formal objects} which are recursively generated by a
process of repetition; the stroke symbols are then concrete representations of
these formal objects (see \cite{Sieg:02}\index{Sieg, W.}).

Sometimes Hilbert's view is presented as if Hilbert claimed that the numbers
are signs on paper.  It is important to stress that this is a
misrepresentation, that the numerals are not physical objects in the sense
that truths of elementary number theory are dependent only on external
physical facts or even physical possibilities (e.g., on what sorts of stroke
symbols it is possible to write down).  \cite{Hilbert:26} made too much of
the fact that for all we know, neither the infinitely small nor the infinitely
large are actualized in physical space and time, yet he certainly held that
the number of strokes in a numeral is at least potentially infinite. It is
also essential to the conception that the numerals are sequences of one kind
of sign, and that they are somehow dependent on being grasped as such a
sequence, that they do not exist independently of our intuition of them.  Only
our seeing or using ``$\mid\mid\mid\mid$'' as a sequence of 4 strokes as
opposed to a sequence of 2 symbols of the form ``$\mid\mid$'' makes
``$\mid\mid\mid\mid$'' into the numeral that it is.  This raises the question
of individuation of stroke symbols.  An alternative account would have
numerals be mental constructions.  However, Bernays\index{Bernays, P.} denied
also this, writing that ``the objects of intuitive number theory, the number
signs, are, according to Hilbert, also not `created by thought'.  But this
does not mean that they exist independently of their {\em intuitive
  construction}, to use the Kantian term that is quite appropriate here''
\cite[226]{Bernays:23b}.  \cite{Kitcher:76}\index{Kitcher, P.} proposes the
view that, whatever the numerals are, the strokes on paper or the stroke
sequences contemplated by the mind merely \emph{represent} the numerals.
According to Hilbert and Bernays\index{Bernays, P.}, the numerals are given in
our representation, but they are not merely subjective ``mental cartoons''
(Kitcher's term).
\begin{quote} 
  If we want [\dots] the ordinal numbers as definite objects free of all
  inessential elements, then in each case we have to take the mere schema of
  the relevant figure of repetition [{\em Wiederholungsfigur}] as an object;
  this requires a very high abstraction.  We are free, however, to represent
  these purely formal objects by concrete objects (``number signs''); these
  contain then inessential, arbitrarily added properties, which, however, are
  also easily grasped as such. \cite[244]{Bernays:30}\index{Bernays, P.}
\end{quote}
One version of this view would be to hold that the numerals are {\em
  types} of stroke-symbols as represented in intuition.  This is
the interpretation that \cite{Tait:81}\index{Tait, W. W.} gives. At first
glance, this seems to be a viable reading of Hilbert.  It takes care of the
difficulties that the reading of numerals-as-tokens (both physical and mental)
faces, and it gives an account of how numerals can be dependent on their
intuitive construction while at the same time not being created by thought.
The reasoning that leads Tait to put forward his reading lies in several
constraints that Hilbert and Bernays\index{Bernays, P.} put on the numerals.
For instance, \cite[159]{Bernays:23b} writes that ``figures [i.e., numerals]
{\em are} not shapes, they {\em have} a shape.''  Thus it is the shape of the
numerals, and not the numerals themselves, which is supposed to be independent
of place and time, independent of the circumstances of production, independent
of inessential differences in execution, and capable of secure recognition in
all circumstances \cite[163]{Hilbert:22a}. Tait\index{Tait, W. W.} infers from
this that identity between numerals is type identity, and hence, that numerals
should be construed as types of stroke symbols.

Types are ordinarily considered to be abstract objects and not located in
space or time. Taking the numerals as intuitive representations of sign types
might commit us to taking these abstract objects as existing independently of
their intuitive representation.  That numerals are ``space- and timeless'' is
a consequence that already \cite[158]{Muller:23}\index{M\"uller, A.} thought
could be drawn from Hilbert's statements, and that was in turn denied by
Bernays\index{Bernays, P.}. The reason is that a view on which numerals are
space- and timeless objects existing independently of us would be committed to
them existing simultaneously as a completed totality, and this is exactly what
Hilbert is objecting to.
\begin{quote}
  It is by no means compatible, however, with Hilbert's basic thoughts to
  introduce the numbers as ideal objects ``with quite different determinations
  from those of sensible objects,'' ``which exist entirely independent of
  us.''  By this we would go beyond the domain of the immediately certain.  In
  particular, this would be evident in the fact that we would consequently
  have to assume the numbers {\em as all existing simultaneously.} But this
  would mean to assume at the outset that which Hilbert considers to be
  problematic.  \cite[225--26]{Bernays:23b}\index{Bernays, P.}
\end{quote}
This is not to say that it is {\em incoherent} to consider the numbers
as being abstract objects, only that the finitary viewpoint prohibits
such a view.  Bernays\index{Bernays, P.} goes on to say:
\begin{quote}
  Hilbert's theory does not exclude the possibility of a philosophical
  attitude which conceives of the numbers [but not the finitist's numerals] as
  existing, non-sensible objects (and thus the same kind of ideal existence
  would then have to be attributed to transfinite numbers as well, and in
  particular to the numbers of the so-called second number class).
  Nevertheless the aim of Hilbert's theory is to make such an attitude
  dispensable for the foundation of the exact sciences.
  \cite[226]{Bernays:23b}
\end{quote}
Another open question in this regard is exactly what Hilbert meant by
``concrete.''  He very likely did not use it in the same sense as it is used
today, i.e., as characteristic of spatio-temporal physical objects in contrast
to ``abstract''\index{abstract objects} objects. However, sign types certainly
are different from full-fledged abstracta like pure sets in that all their
tokens are concrete. Parsons takes account of this difference by using the
term ``quasi-concrete''\index{quasi-concrete} for such abstracta. Tait, on the
other hand, thinks that even the tokens are not concrete physical objects, but
abstract themselves.

Now what is the epistemological status of the finitary objects?  In order to
carry out the task of providing a secure foundation for infinitary
mathematics, access to finitary objects must be immediate and certain.
Hilbert's philosophical background was broadly Kantian, as was
Bernays\index{Bernays, P.}'s, who was closely affiliated with the neo-Kantian
school of philosophy around Leonard Nelson\index{Nelson, L.} in G\"ottingen.
Hilbert's characterization of finitism often refers to Kantian
intuition\index{intuition!Kantian}, and the objects of finitism as objects
given intuitively.  Indeed, in Kant's epistemology, immediacy is a defining
characteristic of intuitive knowledge.  The question is, what kind of
intuition is at play?  \cite{Mancosu:98c}\index{Mancosu, P.} identifies a
shift in this regard.  He argues that whereas the intuition involved in
Hilbert's early papers was a kind of perceptual intuition, in later writings
(e.g., \cite{Bernays:28b}\index{Bernays, P.}) it is identified as a form of
pure intuition\index{intuition!Kantian} in the Kantian sense.  \cite[266--267]{Hilbert:31a} later sees
the finite mode of thought as a separate source of \emph{a priori} knowledge
in addition to pure intuition (e.g., of space) and reason, claiming that he
has ``recognized and characterized the third source of knowledge that
accompanies experience and logic.''  Both Bernays\index{Bernays, P.} and
Hilbert justify finitary knowledge in broadly Kantian terms (without however
going so far as to provide a transcendental deduction), characterizing
finitary reasoning as the kind of reasoning that underlies all mathematical,
and indeed, scientific, thinking, and without which such thought would be
impossible.\footnote{See \cite{Kitcher:76}\index{Kitcher, P.} and
  \cite{Parsons:98}\index{Parsons, C.} for more discussion on the metaphysics
  and epistemology of finitism, and \cite{Sieg:99}\index{Sieg, W.} for
  historical remarks on the development of finitism.}

The simplest finitary propositions are those about equality and inequality of
numerals.  The finite standpoint moreover allows operations on finitary
objects.  Here the most basic is that of concatenation.  The concatenation of
the numerals $\mid\mid$ and $\mid\mid\mid$ is communicated as ``$2 + 3$,'' and
the statement that $\mid\mid$ concatenated with $\mid\mid\mid$ results in the
same numeral as $\mid\mid\mid$ concatenated with $\mid\mid$ by ``$2 + 3 = 3 +
2$.''  In actual proof-theoretic practice, as well as explicitly in
\cite{Bernays:30,HilbertBernays:34}\index{Bernays, P.}, these basic operations
are generalized to operations defined by recursion, paradigmatically,
primitive recursion\index{primitive recursion}, e.g., multiplication and
exponentiation.  Roughly, a primitive recursive definition of a numerical
operation is one in which the function to be defined, $f$, is given by two
equations
\begin{eqnarray*}
f(0, \mathfrak{m}) & = & g(\mathfrak{m}) \\
f(\mathfrak{n}', \mathfrak{m}) & = & h(\mathfrak{n}, \mathfrak{m}, 
f(\mathfrak{n}, \mathfrak{m})),
\end{eqnarray*}
where $g$ and $h$ are functions already defined, and $\mathfrak{n}'$ is the
successor numeral to $\mathfrak{n}$.  Fraktur letters are used here, as in
Hilbert's writings, as meta-variables for numerals.  For instance, if we accept
the function $g(\mathfrak{m}) = \mathfrak{m}$ (the constant function) and
$h(\mathfrak{n}, \mathfrak{m}, \mathfrak{k}) = \mathfrak{m} + \mathfrak{k}$ as
finitary, then the equations above define a finitary function, in this case,
multiplication $f(\mathfrak{n}, \mathfrak{m}) = \mathfrak{n} \times
\mathfrak{m}$.

Similarly, finitary judgments may involve not just equality or inequality but
also basic decidable properties, such as ``is a prime.''  This is finitarily
acceptable as long as the characteristic function of such a property is itself
finitary: For instance, the operation which transforms a numeral to $\mid$ if
it is prime and to $\mid\mid$ otherwise can be defined by primitive recursion
and is hence finitary.  Such finitary propositions may be combined by the
usual logical operations of conjunction, disjunction, negation, but also
bounded quantification.  \cite{Hilbert:26} gives the example of the
proposition that ``there is a prime number between $\mathfrak{p} + 1$ and
$\mathfrak{p}! + 1$'' where $\mathfrak{p}$ is a certain large prime.  This
statement is finitarily acceptable since it ``serves merely to abbreviate the
proposition'' that either $\mathfrak{p} + 1$ or $\mathfrak{p} + 2 $ or
$\mathfrak{p} + 3$ or \dots or $\mathfrak{p}! + 1$ is a prime.

The problematic finitary propositions are those that express general facts
about numerals such as that $1 + \mathfrak{n} = \mathfrak{n} + 1$ for any
given numeral $\mathfrak{n}$.  It is problematic because, as Hilbert puts it,
it ``is from the finitist point of view \emph{incapable of being negated}''
\cite[378]{Hilbert:26}. By this he means that the contradictory proposition
that there is a numeral $\mathfrak{n}$ for which $1 + \mathfrak{n} \neq
\mathfrak{n} + 1$ is not finitarily meaningful.  ``One cannot, after all, try
out all numbers'' \cite[470]{Hilbert:28}.  For the same reason, a finitary
general proposition is not to be understood as an infinite conjunction but
``only as a hypothetical judgment that comes to assert something when a
numeral is given'' \cite[378]{Hilbert:26}. Even though they are problematic in
this sense, general finitary statements are of particular importance to
Hilbert's proof theory, since the statement of consistency of a formal system
$T$ is of such a general form: for any given sequence~$p$ of formulas, $p$ is
not a derivation of a contradiction in~$T$.

Even though in general existential statements are not finitarily meaningful,
they may be \emph{given} finitary meaning if the witness is given by a
finitary function.  For instance, the finitary content of Euclid's theorem
that for every prime $\mathfrak{p}$ there is a prime $> \mathfrak{p}$, is
that given a specific prime $\mathfrak{p}$ one can produce, by a finitary
operation, another prime $> \mathfrak{p}$ (viz., by testing all numbers
between $\mathfrak{p}$ and $\mathfrak{p}! + 1$.  This view is discussed by
\cite{Bernays:30}\index{Bernays, P.} and plays an important role in the uses
\cite{Gentzen:36} and others make of proof theory.

\subsection{Analyses of finitism}

Hilbert's substantial philosophical claims about the finitary standpoint are
difficult to flesh out.  For instance, Hilbert and Bernays\index{Bernays, P.}
both appeal to the role of Kantian intuition for our apprehension of finitary
objects (they are given in the faculty of representation).  Supposing one
accepts this line of epistemic justification in principle, it is plausible
that the simplest examples of finitary objects and propositions, and perhaps
even simple cases of finitary operations such as concatenations of numerals
can be given a satisfactory account.  However, it is unclear exactly which
more complex objects, propositions, and operations should be admitted as
finitary, and how the account can be extended to cover them.  This has led to
substantial debate since the 1920s about the nature of finitary reasoning and
its justification.

Of crucial importance to both an understanding of finitism and of Hilbert's
proof theory is the question of what operations and what principles of proof
should be allowed from the finitist standpoint.  That a general answer is
necessary is clear from the demands of Hilbert's proof theory, i.e., it is not
to be expected that given a formal system of mathematics (or even a single
sequence of formulas) one can ``see'' that it is consistent (or that it cannot
be a genuine derivation of an inconsistency) the way we can see, e.g., that
${\mid\mid} + {\mid\mid\mid} = {\mid\mid\mid} + {\mid\mid}$.  What is required
for a consistency proof\index{consistency proof} is an operation which, given
a formal derivation, transforms such a derivation into one of a special form,
plus proofs that the operation in fact succeeds in every case and that proofs
of the special kind cannot be proofs of an inconsistency.  To count as a
finitary consistency proof, the operation itself must be acceptable from the
finitist standpoint, and the proofs required must use only finitarily
acceptable principles.

Hilbert never gave a general account of which operations and methods of proof
are acceptable from the finitist standpoint, but only examples of operations
and methods of inference in contentual finitary number theory which he
accepted as finitary.  Contentual induction\index{induction!contentual} was
accepted in its application to finitary statements of the hypothetical,
general kind explicitly in \shortcite{Hilbert:22a}.  \cite[1139]{Hilbert:23}
said that intuitive thought ``includes recursion and intuitive induction for
finite existing totalities,'' and used exponentiation in an example in
\citeyear{Hilbert:28}.  \cite{Bernays:30}\index{Bernays, P.} explained how
exponentiation may be understood as a finitary operation on numerals.
\cite{HilbertBernays:34}\index{Bernays, P.} give the only general account of
finitary contentual number theory; according to it, operations defined by
primitive recursion and proofs using induction are finitarily acceptable.  All
of these methods, in their application in the domain of numbers, can be
formalized in a system known as primitive recursive arithmetic\index{primitive
  recursive arithmetic} (\emph{PRA}), which allows definitions of functions by
primitive recursion and induction on quantifier-free formulas.  However,
neither Hilbert nor Bernays\index{Bernays, P.} ever claimed that \emph{only}
primitive recursive\index{primitive recursion!and finitism|(} operations count
as finitary, and non-primitive recursive methods were used in ostensibly
finitary consistency proofs already in 1923 (see \cite{Tait:02}\index{Tait, W.
  W.} and \cite{Zach:03}\index{Zach, R.}).  These include, in particular, the
consistency proof\index{consistency proof} of a formal system of number theory
corresponding to primitive recursive arithmetic in \cite{Hilbert:22b}, and a
stronger system in Ackermann\index{Ackermann, W.}'s dissertation
\cite{Ackermann:24}.\footnote{Ackermann in fact used transfinite
  induction\index{induction!transfinite} up to $\omega^{\omega^\omega}$ to
  justify a non-primitive recursive reduction procedure.}

Although Hilbert and his collaborators used methods which go beyond the
primitive recursive and accepted them as finitary, it is still unclear whether
they (a) realized that these methods were not primitive recursive and (b)
whether they would still have accepted them as finitary if they
had.\footnote{See Tait's discussion in the Appendix to Chapters~1 and 2 in
  \cite{Tait:05}.}  The conceptual issue is which operations \emph{should} be
considered as finitary.  Since Hilbert was less than completely clear on what
the finitary standpoint consists in, there is some leeway in setting up the
constraints, epistemological and otherwise, an analysis of finitist operation
and proof must fulfill.  Hilbert characterized the objects of finitary number
theory as ``intuitively given,'' as ``surveyable in all their parts,'' and
said that their having basic properties must ``exist intuitively'' for us.
\cite[216]{Bernays:22}\index{Bernays, P.} suggests that in finitary
mathematics, only ``primitive intuitive cognitions come into play,'' and uses
the term ``point of view of intuitive evidence'' in connection with finitism
\cite[250]{Bernays:30}\index{Bernays, P.}. This characterization of finitism
as primarily to do with \emph{intuition}\index{intuition!and finitism} and
intuitive knowledge has been emphasized in particular by
\cite{Parsons:98}\index{Parsons, C.} who argues that what can count as
finitary on this understanding is not more than those arithmetical operations
that can be defined from addition and multiplication using bounded recursion.
In particular, according to Parsons, exponentiation and general primitive
recursion are not finitarily acceptable.

The thesis that finitism coincides with primitive recursive
reasoning\index{primitive recursion} has
received a forceful and widely accepted defense by \cite{Tait:81}.
Tait\index{Tait, W. W.|(}, in contrast to Parsons\index{Parsons, C.}, rejects
the aspect of representability in intuition as the hallmark of the finitary;
instead he takes finitary reasoning to be ``a minimal kind of reasoning
supposed by all non-trivial mathematical reasoning about numbers'' and
analyzes finitary operations and methods of proof as those that are implicit
in the very notion of number as the form of a finite sequence.  This analysis
of finitism is supported by Hilbert's contention that finitary reasoning is a
precondition for logical and mathematical, indeed, any scientific thinking
\cite[267]{Hilbert:31a}.\index{primitive recursion!and finitism|)} The crucial
difference between Tait's conception of finitism and Parsons (as well as
Hilbert's own) is that according to Tait there is no ultimate epistemological
foundation for finitism which yields the security of finitary reasoning for
which Hilbert appealed to intuition.  Tait argues that
\begin{quote}
[\dots] no absolute conception of security is realized by finitism
or any other kind of mathematical reasoning.  Rather, the special role
of finitism consists in the circumstance that it is a minimal kind of
reasoning presupposed by all nontrivial mathematical reasoning about
numbers.  And for this reason it is \emph{indubitable} in a Cartesian
sense that there is no preferred or even equally preferable ground on
which to stand and criticize it. Thus finitism is fundamental to
number-theoretical mathematics even if it is not a foundation in the
sense Hilbert wished. \cite[525]{Tait:81}
\end{quote}

Another interesting analysis of finitary proof, which, however, does not
provide as detailed a philosophical justification, was proposed by
\cite{Kreisel:58}\index{Kreisel, G.|(}.  It yields the result that exactly
those functions are finitary which can be proved to be well-defined in
first-order arithmetic \emph{PA}.\footnote{\cite[Section
  3.5]{Kreisel:70}\index{Kreisel, G.} provides another analysis by focusing on
  what is ``visualizable.'' The result is the same: finitary functions turn
  out to be just those provably total in~\emph{PA}.}  Kreisel's proposal makes
use of the notions of formalizations of provability predicates and ordinal
progressions of theories.  Kreisel argues that if $\mathit{Pr}(\Godelnum{A})$
has been recognized to be a provability predicate for a partial formalization
$\Sigma_\mu$ of finitary reasoning, and
$\mathit{Pr}(\Godelnum{A(0^{(x)})})$\footnote{Here, $x$ is a free variable,
  and $\Godelnum{A(0^{(x)})}$ is the function of $x$ which computes
  $\Godelnum{A(0^{\prime\cdots\prime})}$ with $x$ occurrences of~$'$.} is
provable in $\Sigma_\mu$ (and hence established by finitary means), then the
finitist is entitled to also accept $A(x)$ as finitarily established.  If that
is the case, we may add $A(x)$ as an axiom to $\Sigma_\mu$ to obtain a new
theory $\Sigma_\nu$ which is also finitarily justified.  On Kreisel's account,
finitary provability coincides with the provability in some $\Sigma_\nu$ so
obtained, starting from $\Sigma_0 = \mathit{PRA}$\index{primitive recursive
  arithmetic}.  If some $\Sigma_\nu$ proves $\exists y\, A(x, y)$, for $A(x,
y)$ a primitive recursive equation, then $f(x) =$ the least $y$ such that
$A(x, y)$ is finitarily justified.  Kreisel sketches a proof of the theorem
that the functions so justified are exactly those which are provably total in
$\mathit{PA}$, and hence there are finitary functions which are not primitive
recursive.

\cite[\S13]{Tait:81} also contains a discussion of Kreisel's analysis.  In
order to obtain a non-primitive recursive function on Kreisel's account, we
must properly extend~$\Sigma_0$ since the provably total functions of
$\Sigma_0 = \mathit{PRA}$\index{primitive recursive arithmetic} are just the
primitive recursive functions.  So suppose we have that
$\mathit{PRA}$\index{primitive recursive arithmetic} proves the
arithmetization of the claim that $\exists y\, A(0^{(x)}, y)$ is
provable.\footnote{In other words, there are primitive recursive functions
  $h(x)$ and $g(x)$ so that $\mathit{PRA}$ proves that $g(x)$ codes a
  derivation in $\mathit{PRA}$ of the formula $\Godelnum{A(0^{(x)}, t)}$,
  where $t$ is the primitive recursive term (containing only the free variable
  $x$) which is coded by $h(x)$.} This, according to Kreisel, justifies the
acceptance of $f(x)$ as defined above as finitary, because a finitary proof of
the general fact that it is provable in $\mathit{PRA}$ that $f(x)$ is defined,
together with the acceptance of $\mathit{PRA}$\index{primitive recursive
  arithmetic} as finitarily acceptable, amounts to a finitary proof that
$f(x)$ is defined for all~$x$.  The weak point in this justification,
according to Tait, is this:
\begin{quote}
  For the finitist to recognize the validity of primitive recursive
  arithmetic, he must recognize the \emph{general} validity of definition of
  functions by primitive recursion.  But he cannot even formulate this since
  it involves the notion of function.
\end{quote}
Tait's point here is that there is a significant difference between accepting
each primitive recursive definition individually as finitary, and accepting
primitive recursion in general as a finitarily valid principle.  The finitist
is able to do the former, but not the latter.  For to accept primitive
recursion in general as a finitarily valid principle of definition, one would
either, as Tait puts it, need to invoke the notion of a function (which is not
a finitary object), or one needs a justification for why, say, all primitive
recursive \emph{terms} are calculable for every argument---and for this a
finitary evaluation procedure for primitive recursive terms is necessary.
Such an evaluation procedure, however, cannot be primitive recursive.  And
prior to the extension of $\mathit{PRA}$\index{primitive recursive arithmetic}
to include the new non-primitive recursive function $f(x)$ there is no reason
to suppose that such an evaluation procedure exists.  Although Tait's
objection is directed at Kreisel's analysis of finitary function, it of course
also raises doubts about Kreisel's account of finitary proof.\index{Kreisel,
  G.|)}\index{finitism|)}\index{Tait, W. W.|)}

\subsection{Hilbert's program and instrumentalism}

Hilbert's program has often been interpreted as an
instrumentalist\index{instrumentalism!and Hilbert's program|(}\index{Hilbert's
  program!and instrumentalism|(} account of mathematics.  This reading relies
on the distinction Hilbert makes between the finitary part of mathematics and
the non-finitary rest which is in need of grounding (via finitary
meta-mathematics).  The finitary part Hilbert calls ``contentual,'' i.e., its
propositions and proofs have content.  The infinitary part, on the other hand,
is ``not meaningful from a finitary point of view.''  This distinction
corresponds to a distinction between formulas of the the axiomatic systems of
mathematics for which consistency proofs\index{consistency proof} are being
sought.  Some of the formulas correspond to contentual, finitary propositions:
they are the ``real'' formulas.  The rest are called ``ideal.''  They are
added to the real part of our mathematical theories in order to preserve
classical inferences such as the principle of the excluded middle for infinite
totalities, i.e., the principle that either all numbers have a given property
or there is a number which does not have it.  This disjunction is not
finitarily valid, as we saw above.  Hilbert first mentioned ``ideal''
propositions in \shortcite{Hilbert:26}, although the distinction was hinted at
in~\shortcite{Hilbert:23}.  In the latter paper, Hilbert presented a formal
system of quantifier-free number theory about which he says that ``the
provable formulae we acquire in this way all have the character of the
finite'' (1139).  Then the transfinite axioms (i.e., quantifiers) are added to
simplify and complete the theory (1144). Here he draws the analogy with the
method of ideal elements: ``In my proof theory, the transfinite axioms and
formulae are adjoined to the finite axioms, just as in the theory of complex
variables the imaginary elements are adjoined to the real, and just as in
geometry the ideal constructions are adjoined to the actual'' (ibid). When
Hilbert, in \shortcite{Hilbert:26}, explicitly introduces the notion of an
ideal proposition, and in \shortcite{Hilbert:28}, when he first speaks of
\emph{real propositions} in addition to the ideal, he is quite clear that the
real part of the theory consists only of decidable, variable-free formulas.
They are supposed to be ``directly capable of verification''---akin to
propositions derived from laws of nature which can be checked by experiment
\cite[475]{Hilbert:28}.\footnote{This reading is not universally accepted.
  \cite{Detlefsen:90}, for instance, considers the real formulas to also
  include the general formulas, i.e., formulas with free variables.  See
  \cite{Zach:02}\index{Zach, R.} for a defense of the reading given here.}  It
is this extension of the real part of the theory by the ideal, infinitary part
that is in need of justification by a consistency proof\index{consistency
  proof}:
\begin{quote}
  For there is a condition, a single but absolutely necessary one, to which
  the use of the method of ideal elements is subject, and that is the
  \emph{proof of consistency}; for, extension by the addition of ideals is
  legitimate only if no contradiction is thereby brought about in the old,
  narrower domain, that is, if the relations that result for the old objects
  whenever the ideal objects are eliminated are valid in the old domain.
  \cite[383]{Hilbert:26}
\end{quote}
\cite{Weyl:25}\index{Weyl, H.} described Hilbert's project as replacing
meaningful mathematics by a meaningless game of formulas.  He noted that
Hilbert wanted to ``secure not \emph{truth}, but the \emph{consistency} of
analysis'' and suggested a criticism that echoes an earlier one by
Frege\index{Frege, G.}: Why should we take consistency of a formal system of
mathematics as a reason to believe in the truth of the pre-formal mathematics
it codifies?  Is Hilbert's meaningless inventory of formulas not just ``the
bloodless ghost of analysis''?  Weyl\index{Weyl, H.} suggested a solution:
\begin{quote}
  [I]f mathematics is to remain a serious cultural concern, then some
  \emph{sense} must be attached to Hilbert's game of formulae, and I see only
  one possibility of attributing to it (including its transfinite components)
  an independent intellectual meaning. In theoretical physics we have before
  us the great example of a [kind of] knowledge of completely different
  character than the common or phenomenal knowledge that expresses purely what
  is given in intuition. While in this case every judgment has its own sense
  that is completely realizable within intuition, this is by no means the case
  for the statements of theoretical physics.  In that case it is rather
  \emph{the system as a whole} that is in question if confronted with
  experience.  \cite[140]{Weyl:25}\index{Weyl, H.}
\end{quote}  
The analogy with physics is striking. Hilbert himself used a similar analogy
in \shortcite{Hilbert:28}. He there suggested that consistency is not the
only virtue ideal mathematics has: transfinite inference simplifies and
abbreviates proofs, ``brevity and economy of thought are the \emph{raison
  d'\^etre} of existence proofs'' (476).  The formal system of transfinite
logic is not arbitrary: ``This formula game is carried out according to
certain definite rules, in which the \emph{technique of our thinking} is
expressed. [...] The fundamental idea of my proof theory is none other than to
describe the activity of our understanding, to make a protocol of the rules
according to which our thinking actually proceeds'' (ibid).

Although these remarks are suggestive, they do not force an interpretation of
Hilbert as an instrumentalist.  Most commentators have taken this reading
(including \cite{Kitcher:76}\index{Kitcher, P.},
\cite{Resnik:80}\index{Resnik, M.}, \cite{Giaquinto:83}\index{Giaquinto, M.},
\cite{Sieg:90}\index{Sieg, W.}, and in particular \cite{Detlefsen:86}) in that
they interpret Hilbert's explanation that the ideal propositions ``have no
meaning in themselves'' \cite[381]{Hilbert:26} as claiming that classical
mathematics is a \emph{mere} instrument, and that statements of transfinite
mathematics have no truth value.  To the extent that this is accurate,
however, it must be understood as a methodological instrumentalism: A
successful execution of the proof-theoretic program would show that one could
pretend \emph{as if} mathematics was meaningless.\footnote{On this point see
  also \cite{Sieg:99}\index{Sieg, W.}, esp.~B3 and the conclusion, and
  \cite{Sieg:02}.}  The analogy with physics is therefore not: transfinite
propositions have no meaning just as propositions involving theoretical terms
have no meaning, but: transfinite propositions require no direct intuitive
meaning just as one does not have to directly see electrons in order to
theorize about them.  \cite{Hallett:90}\index{Hallett, M.}, taking into
account the 19th century mathematical background from which Hilbert came as
well as published and unpublished sources from Hilbert's entire career (in
particular \cite{Hilbert:92}, the most extensive discussion of the method of
ideal elements) comes to the following conclusion:
\begin{quote}
  [Hilbert's treatment of philosophical questions] is \emph{not} meant as a
  kind of instrumentalist agnosticism about existence and truth and so forth.
  On the contrary, it is meant to provide a non-skeptical and positive
  solution to such problems, a solution couched in cognitively accessible
  terms.  And, it appears, the same solution holds for both mathematical and
  physical theories.  Once new concepts or ``ideal elements'' or new
  theoretical terms have been accepted, then they exist in the sense in which
  \emph{any} theoretical entities exist. \cite[239]{Hallett:90}\index{Hallett,
    M.}
\end{quote}
This conclusion is in line with Weyl\index{Weyl, H.}'s assessment in
\shortcite{Weyl:28}. When Weyl eventually turned away from
intuitionism\index{intuitionism},\footnote{For the reasons for Weyl's
  rejection of intuitionism, see \cite{Mancosu:02b}\index{Mancosu,
    P.}\index{Ryckman, T.}.} he
emphasized the purpose of Hilbert's proof theory, not to turn mathematics into
a meaningless game of symbols, but to turn it into a theoretical science which
codifies scientific (mathematical) practice.

The reading of Hilbert as an instrumentalist goes hand in hand with a reading
of the proof-theoretic program as a reductionist project.  The instrumentalist
reading interprets ideal mathematics as a meaningless formalism, which
simplifies and ``rounds out'' mathematical reasoning.  But a consistency
proof\index{consistency proof} of ideal mathematics by itself does not explain
what ideal mathematics is an instrument \emph{for}.  Thus, commentators have
sought to elaborate on Hilbert's discussion of consistency by pointing out
that consistency proofs\index{consistency proof} do not just establish that
ideal theories are free from formal contradictions, but that they establish
more than mere consistency.  They establish
\emph{conservativity}\index{conservativity|(} of the ideal over the real part
of the theory, in the following sense: If the ideal theory proves a real
statement, then the real statement is also provable by real, finitary means.
Such reductivist projects were common in the philosophy of science at the
time, as was pointed out by \cite{Giaquinto:83}\index{Giaquinto, M.}.  A
conservativity proof justifies the use of transfinite mathematics: it is not
only internally consistent, but it proves only true intuitive propositions.

On this picture, classical mathematics is to be formalized in a system which
includes formalizations of all the directly verifiable (by calculation)
propositions of contentual finite number theory.  The consistency
proof\index{consistency proof} should show that all real propositions which
can be proved by ideal methods are true, i.e., can be directly verified by
finite calculation.  Actual proofs such as the
$\epsilon$-substitution\index{epsilon substitution@$\epsilon$-substitution
  method} procedure are of such a kind: they provide finitary procedures which
eliminate transfinite elements from proofs of real statements. In particular,
they turn putative ideal derivations of $0 = 1$ into derivations in the real
part of the theory; the impossibility of such a derivation establishes
consistency of the theory. Indeed, Hilbert saw that something stronger is
true: not only does a consistency proof\index{consistency proof} establish
truth of real formulas provable by ideal methods, but it yields finitary
proofs of finitary \emph{general} propositions if the corresponding
free-variable formula is derivable by ideal methods \cite[474]{Hilbert:28}.

It bears pointing out that Hilbert never clearly articulated conservativity of
the ideal over the real for finitary general statements as an aim of his
foundational project.  There are contemporary commentators, e.g.,
\cite{Neumann:31}\index{von Neumann, J.} who attribute to Hilbert an interest
in conservativity proofs, but it was only Bernays\index{Bernays, P.} in the
\emph{Grundlagen der Mathematik} who pointed out that consistency
proofs\index{consistency proof} themselves established not only the truth of
variable-free formulas provable by ideal methods, but also of free-variable
theorems.  In this context, Bernays\index{Bernays, P.} used the term
`verifiable' (\emph{verifizierbar}): a free-variable formula $A(x)$ is
verifiable if each numerical instance $A(\mathfrak{z})$ is (finitarily) true.
The proof transformation methods used in consistency proofs\index{consistency
  proof} for systems of arithmetic in \cite[248,
298]{HilbertBernays:34}\index{Bernays, P.} can be applied not only to putative
proofs of $0 = 1$, but generally to proofs of formulas with free-variables.
If we have a proof of $A(x)$, then the following method constitutes a finitary
proof that, for any $\mathfrak{z}$, $A(\mathfrak{z})$ is true.  From the
derivation of $A(x)$ we obtain a derivation of $A(\mathfrak{z})$ by
substitution. The procedure given in the consistency proof\index{consistency
  proof} transforms this derivation into a variable-free derivation of
$A(\frak{z})$ in the real part of the theory, which codifies a finitary
calculation that $A(\mathfrak{z})$ is true.

Kreisel\index{Kreisel, G.} was most influential in promoting the
interpretation of the aim of Hilbert's program as an attempt to establish
conservativity of the ideal theory for finitary general propositions all
along.  \cite{Kreisel:51} cites Bernays\index{Bernays, P.}'s results; but in
\cite{Kreisel:58} and later, he instead points to an argument in
\cite[474]{Hilbert:28}. This argument, unlike Bernays\index{Bernays, P.}'s,
does not rely on a particular form of the consistency proof\index{consistency
  proof}.  It assumes only that a finitary consistency proof for an ideal
theory is available.  Assume there is a derivation of $A(x)$.  Now suppose
that for a given $\mathfrak{z}$, $A(\frak{z})$ is not true.  Then $\neg
A(\mathfrak{z})$ would be true,\footnote{This inference uses tertium non
  datur, but only regarding the unproblematic finitary
  statement~$A(\mathfrak{z})$.}  and so there would be a derivation of $\neg
A(\mathfrak{z})$ in the ideal theory (which includes all real theorems).  But
from the derivation of $A(x)$ we obtain, by substitution, a derivation of
$A(\mathfrak{z})$, and hence a contradiction. Since we assume that we have a
finitary consistency proof of~$T$, this cannot be the case.  Hence,
$A(\mathfrak{z})$ must be true.\index{instrumentalism!and Hilbert's
  program|)}\index{Hilbert's program!and
  instrumentalism|)}\index{conservativity|)}

\subsection{Hilbert's program and G\"odel's incompleteness theorems}

G\"odel\index{G\"odel, K.} announced the second incompleteness
theorem\index{incompleteness theorems|(}\index{incompleteness theorems!second
  incompleteness theorem} in an abstract published in October 1930: no
consistency proof\index{consistency proof} of systems such as
\emph{Principia}\index{Principia Mathematica}, Zermelo-Fraenkel set
theory\index{set theory}, or the systems investigated by
Ackermann\index{Ackermann, W.} and von Neumann\index{von Neumann, J.} is
possible by methods which can be formulated in these systems.  In the full
version of his paper, \cite{Godel:31} left open the possibility that there
could be finitary methods which are not formalizable in these systems and
which would yield the required consistency proofs\index{consistency proof}.
Bernays\index{Bernays, P.}'s first reaction in a letter to
G\"odel\index{G\"odel, K.} in January 1931 was likewise that ``if, as von
Neumann\index{von Neumann, J.} does, one takes it as certain that any and
every finitary consideration may be formalized within the system $P$---like
you, I regard that in no way as settled---one comes to the conclusion that a
finitary demonstration of the consistency of $P$ is impossible''
\cite[87]{Godel:03}.

Through a careful (``G\"odel''-) coding of sequences of symbols (formulas,
proofs) by numbers, G\"odel\index{G\"odel, K.} showed that in theories $T$
which contain a sufficient amount of arithmetic, it is possible to produce a
formula $\mathit{Pr}(x, y)$ which expresses that $x$ is (the code of) a proof
of (the formula with code)~$y$.  Specifically, if $\Godelnum{0 = 1}$ is the
code of the formula $0 = 1$, then $\mathit{Con}_T \equiv \forall x\, \neg
\mathit{Pr}(x, \Godelnum{0 = 1})$ expresses that $T$ is consistent (no number
is the code of a derivation in $T$ of $0 = 1$). The second incompleteness
theorem\index{incompleteness theorems!second incompleteness theorem} (G2) says
that under certain assumptions about $T$ and the coding apparatus, $T$ does
not prove $\mathit{Con}_T$.  Now suppose there were a finitary consistency
proof\index{consistency proof} of $T$. The methods used in such a proof would
presumably be formalizable in~$T$.  (``Formalizable'' means that, roughly, if
the proof uses a finitary operation $f$ on derivations which transforms any
derivation~$d$ into a derivation $f(d)$ of a simple form; then there is a
formula~$F(x, y)$ so that, for all derivations $d$, $T \vdash F(\Godelnum{d},
\Godelnum{f(d)})$.)  The consistency of~$T$ would be finitarily expressed as
the general hypothetical claim that, if $d$ is any given sequence of symbols,
$d$ is not a derivation in $T$ of the formula~$0 = 1$.  The formalization of
this proposition is the formula $\neg \mathit{Pr}(x, \ulcorner 0 =
1\urcorner)$ in which the variable~$x$ occurs free.  If there were a finitary
proof of the consistency of~$T$, its formalization would yield a derivation in
$T$ of $\neg \mathit{Pr}(x, \Godelnum{0 = 1})$, from which $Con_T$ can be
derived in $T$ by simple universal generalization on~$x$.  Yet, a derivation
of $\mathit{Con}_T$ in $T$ is ruled out by G2.

G\"odel\index{G\"odel, K.} and Bernays\index{Bernays, P.} initially thought
that the difficulty for the consistency proof of Peano arithmetic
$\mathit{PA}$\index{Peano arithmetic} could be overcome by employing methods
which, although not formalizable in $\mathit{PA}$, are nonetheless finitary.
Bernays\index{Bernays, P.} did not seem to have any particular candidates for
such a method in mind, and he thought that all methods which were up to then
employed in finitary considerations were in fact formalizable in
$\mathit{PA}$.  Another option he considered was an extension of the notion of
an axiomatic theory by a finitary version of the $\omega$-rule\index{omega
  rule@$\omega$-rule} proposed by Hilbert \shortcite{Hilbert:31b,Hilbert:31a}.
Such theories might avoid the impact of G\"odel\index{G\"odel, K.}'s
incompleteness theorem since they do not satisfy the conditions of the
incompleteness theorems: the set of axioms would not be decidable.  It is fair
to say, however, that since about 1934 it has been almost universally accepted
that the methods of proof accepted as finitary prior to G\"odel\index{G\"odel,
  K.}'s results are all formalizable in \emph{PA} and that the incompleteness
theorems do show that there can be no finitary consistency
proofs\index{consistency proof} for axiomatic theories of the kind originally
considered by Hilbert.

The reaction to the incompleteness theorems in the Hilbert school then focused
on extensions of the original finitary standpoint in which consistency
proofs\index{consistency proof} for substantial theories like $\mathit{PA}$
can be carried out. Such extensions have been proposed and defended on broadly
finitary grounds, e.g., \cite{Gentzen:36} defended the use of transfinite
induction up to $\epsilon_0$ in his consistency proof for \emph{PA} as
``indisputable,'' and \cite{Takeuti:87}\index{Takeuti, G.} gave another
defense.  In the Gentzen\index{Gentzen, G.}-Sch\"utte tradition of proof
theory by ordinal analysis\index{ordinal analysis}, the proof methods used to
give consistency proofs\index{consistency proof} are all of this sort. To wit,
one uses transfinite induction on ordinal notation systems which name larger
and larger ordinals.  The more complicated the ordering, the more difficult it
is to see that the induction principle in question is finitarily justified.
Another extension of the finitary standpoint is due to \cite{Godel:58}.

\cite{Smorynski:77}\index{Smory\'nski, C.}, following earlier suggestions by
Kreisel\index{Kreisel, G.}, has argued that already the first incompleteness
theorem\index{incompleteness theorems!first incompleteness theorem} defeats
Hilbert's program.  This argument uses the fact that a finitary consistency
proof\index{consistency proof} of an ideal theory~$T$ implies the
conservativity of~$T$ over finitary, real mathematics for general finitary
statements of the form $A(x)$ (with free variable~$x$).  Now
G\"odel\index{G\"odel, K.}'s first incompleteness theorem (G1) states that for
any sufficiently strong, consistent formal theory~$S$ there is a
sentence~$G_S$ which is not derivable in~$S$ if $S$ is consistent.  $G_S$ is a
general real sentence.  Consider a theory $T$ which formalizes ideal
mathematics and contains the theory $S$, which formalizes real mathematics, as
a subtheory.  $S$ satisfies the conditions of G1 and hence $S$ does not
derive~$G_S$.  Yet~$T$, being a formalization of all of mathematics, proves
(via a formalization of G1) that if $S$ is consistent, then $G_S$, but it also
proves the consistency (indeed, the soundness) of $S$.  Hence $T$ proves
$G_S$.  Thus, we have a true real statement which is provable in ideal
mathematics but not in real mathematics.\footnote{The argument appeals to a
  number of perhaps contentious assumptions, such as that $T$ proves the
  soundness of $S$.  For a critique, see~\cite{Detlefsen:90}\index{Detlefsen,
    M.}.}\index{incompleteness theorems|)}\index{Hilbert, D|)}

\section{Hilbert's program now}

\subsection{Detlefsen's Hilbertian instrumentalism}

Detlefsen \shortcite{Detlefsen:79,Detlefsen:86,Detlefsen:01}\index{Detlefsen,
  M.|(}\index{instrumentalism!and Hilbert's program|(}\index{Hilbert's
  program!and instrumentalism|(} has proposed a wide-ranging instrumentalist
account of mathematics based on Hilbert's program which is designed to escape
the difficulties that G\"odel\index{G\"odel, K.}'s incompleteness
theorems\index{incompleteness theorems} poses for Hilbert's original
consistency project.  The project has several parts: \cite{Detlefsen:86} first
gives a detailed analysis and defense of a general instrumentalist view of
mathematics along Hilbertian lines.  This includes an analysis of the
epistemic import of ideal proofs of real statements, which answers a question
that was hardly addressed by Hilbert\index{Hilbert, D.}, either in his mature
writings in the 1920s or in his exchange with Frege\index{Frege, G.} on
formalism and consistency.  This is the question of how manipulation of
meaningless strings of symbols can ever lead to knowledge (of finitary
truths).  Detlefsen then analyzes the characteristics of formal systems of
ideal mathematics \emph{qua} instruments.  On Detlefsen's view, even though,
say, full set theory is commonly accepted as a formalization of infinitary
mathematics, only parts of set theory are in fact instrumentally useful.  In
particular, (1) ideal proofs of real theorems which are more complex than any
real proof of the same theorem do not yield an instrumental advantage, and
hence are not instrumentally useful; and (2) ideal proofs which are too long
or complex to be comprehended by humans, and hence never play a role in actual
mathematical reasoning, are also of no instrumental value.  A proof theoretic
justification of instrumental mathematics, i.e., the proof of the
conservativity of the ideal theory over real mathematics, is only required, so
Detlefsen, for the instrumentally useful part.  Detlefsen introduces the term
``Hilbertian residue'' for that part of ideal mathematics that is
instrumentally useful and hence in need of proof-theoretic justification.  On
his view, then, we only need a consistency proof\index{consistency proof} for
the Hilbertian residue, not for all of ideal mathematics.

This move from a required justification for all of infinitary
mathematics to a justification of only the Hilbertian residue is one step
toward Detlefsen's defense of instrumentalism against G\"odel\index{G\"odel,
  K.}'s incompleteness theorems\index{incompleteness theorems}.  For the
incompleteness theorems only apply under certain conditions, e.g., only when
the theory in question contains enough basic arithmetic to carry out
G\"odel\index{G\"odel, K.} coding, formalization of a proof predicate, and to
prove the diagonal lemma.  The Hilbertian residue of a theory, however, need
not contain a sufficiently strong arithmetical subtheory because of (1) above.
This provides part of Detlefsen's defense against the challenge posed by the
first incompleteness theorem\index{incompleteness theorems!first
  incompleteness theorem} \cite[Appendix]{Detlefsen:86}.  \cite{Detlefsen:90}
also argues that instrumentalism escapes the argument from G1 by denying that
ideal mathematics must be conservative over the real part.  According to him,
Hilbertian instrumentalism requires only that the ideal theory not prove
anything which is in conflict with the real theory; it is not required that
all its real theorems are also provable by real means.  If this defense is
successful, Detlefsen is right to claim that not G1, but only G2 poses a real
challenge to instrumentalism.

Detlefsen presents several lines of defense against the argument from G2, one
of which \shortcite{Detlefsen:79} echoes \cite{Hilbert:31a}\index{Hilbert,
  D.}.  If a version of the $\omega$-rule\index{omega rule@$\omega$-rule} is
finitarily acceptable, then we would have found a finitarily acceptable method
of proof which is not capable of formalization.  Hence, real mathematics is
not a subtheory of the ideal instrument, but this was an assumption necessary
to draw the conclusion that there can be no real consistency
proof\index{consistency proof} of the ideal theory.
\cite{Ignjatovic:94}\index{Ignjatovi\v{c}, A.} raised serious doubts about the
acceptability of Detlefsen's version of the $\omega$-rule, however.
Detlefsen's other argument against the common interpretation of
G\"odel\index{G\"odel, K.}'s second theorem focuses on the notion of
formalization. That the particular formalization of ``$T$ is consistent'' by
G\"odel\index{G\"odel, K.}'s formula $\mathit{Con}_T$ is not provable does not
imply that there could not be other formulas, which \emph{are} provable
in~$T$, and which have as much right to be called ``formalizations of the
consistency of~$T$.''  These rely on different formalizations of the
provability predicate $\mathit{Pr}_T$ than the standard ones.  It is known
that formalized consistency statements are unprovable whenever the provability
predicate obeys certain general derivability conditions.  Detlefsen argues
that these conditions are not necessary for a predicate to count as a genuine
provability predicate, and indeed there are provability predicates which
violate the provability conditions and which give rise to consistency formulas
which \emph{are} provable in their corresponding theories.  These, however,
depend on non-standard conceptions of provability which would likely not have
been accepted by Hilbert\index{Hilbert, D.}.  One quite basic example is the
use of Rosser provability instead of ordinary provability.  On this approach,
a derivation of a formula $A$ only counts as a proof if no derivation with
smaller G\"odel number is a derivation of $\lnot A$.  If
$\mathit{Prov}(x, \Godelnum{A})$ is the standard formalization of ``$x$ is the
code of a derivation of the formula~$A$,'' then the Rosser
provability\index{Rosser provability} predicate is given by
\[
\mathit{RPr}(\Godelnum{A}) \equiv \exists x(\mathit{Prov}(x, \Godelnum{A})
  \land \forall y<x\, \lnot \mathit{Prov}(y, \Godelnum{\lnot A})).
\] 
For this provability predicate, $\lnot \mathit{RPr}(\Godelnum{0=1})$ \emph{is}
provable in, e.g., first-order Peano arithmetic\index{Peano arithmetic}.
Provability of a formula $A$, however, is no longer just a matter of deriving
it from the axioms; one also has to check that all shorter derivations do not
end in~$\lnot A$. Other ``consistency minded'' theories\index{consistency
  minded theory} which prove their own consistency are discussed, e.g., in
\cite{Jeroslow:71,Jeroslow:75}\index{Jeroslow, R.} and especially
\cite{Visser:89}\index{Visser, A.}. The Rosser provability predicate is
studied in, e.g., \cite{Guaspari:79}\index{Guaspari, D.}\index{Solovay, R. M.}
and \cite{Arai:90}\index{Arai, T.}.\footnote{For technical background,
  discussion of intensional provability predicates and examples, see
  \cite{Feferman:60}\index{Feferman, S.}. For discussion, see also
  \cite{Resnik:74a}, \cite{Auerbach:85,Auerbach:92}\index{Auerbach, D.} and
  \cite{Steiner:91}.}

Another interesting aspect of Detlefsen's approach to instrumentalism and
Hilbert's program related to technical work in proof theory is the emphasis on
instrumental utility of ideal proofs.  Hilbert\index{Hilbert, D.}, as we saw
above, himself noted the theoretical and cognitive advantage of ideal methods,
such as increased simplicity of proofs.  In Detlefsen's instrumentalism, such
considerations take center stage.  Even if it is conceded that
G\"odel\index{G\"odel, K.}'s theorems call the success of instrumentalism in
its most general form into question, it would still be of substantial interest
to study restricted cases of conservative extensions of real mathematics which
are instrumentally useful.  To flesh out the notion of ``instrumental
usefulness,'' one obvious characteristic of formal proofs is length.  For
instance, one might take an ideal theory to be useful if its proofs are
substantially shorter than proofs in, say, $\mathit{PRA}$\index{primitive
  recursive arithmetic} of the same theorems.  This question is amenable to
precise proof theoretical study. \cite{Ignjatovic:03}\index{Ignjatovi\v{c},
  A.}\index{Caldon, P.}  prove some related, but on the whole, negative
results: The subsystem of first-order arithmetic $I\Sigma_1$ in which
induction is limited to $\Sigma_1$ formulas has super-exponential
``speed-up''\index{speed up}\index{proof theory!speed up} over
\textit{PRA}\index{primitive recursive arithmetic}. This indicates that using
induction over non-finitary formulas ($\Sigma_1$ formulas have unbounded
existential quantifiers) yields significantly shorter proofs than proofs
without.  However, more comprehensive theories ($\mathit{RCA}$,
$\mathit{WKL}$, see below) which contain some second-order apparatus, do not
significantly shorten proofs vis-\`{a}-vis~$I\Sigma_1$.\index{Detlefsen,
  M.|)}\index{instrumentalism!and Hilbert's program|)}\index{Hilbert's
  program!and instrumentalism|)}

\subsection{Generalized Hilbert programs}

The work of Gentzen\index{Gentzen, G.} on consistency proofs\index{consistency
  proof} for mathematical theories using methods that go beyond the strictly
finitary marks the beginning of an important line of proof-theoretic research.
As outlined in \ref{impact-godel} above, Gentzen\index{Gentzen, G.}'s approach
was to retain the aim of Hilbert's program, viz., to give consistency
proofs\index{consistency proof} for strong mathematical theories by restricted
means.  Because of G\"odel\index{G\"odel, K.}'s incompleteness
theorems\index{incompleteness theorems}, these restricted means are
necessarily not themselves formalizable in the theories whose consistency is
established by them.  Nevertheless, they retain a constructive character, and
attempts have been made to justify them on finitary grounds.

The consistency proof\index{consistency proof} of \cite{Gentzen:36}, as
discussed above, uses the principle of transfinite
induction\index{induction!transfinite} up to~$\epsilon_0$\index{epsilon
  zero@$\epsilon_0$} in order to establish the consistency of first-order
Peano arithmetic\index{Peano arithmetic}.  Gentzen\index{Gentzen, G.}'s use of
a system of notations for ordinals less than~$\epsilon_0$, and the proof of
the termination of a reduction procedure for derivations in~$\mathit{PA}$
based on induction on these ordinal notations, provide the model for the proof
theoretic analysis of axiomatic systems along these lines.  In order to give
an ``ordinal analysis''\index{ordinal analysis}\index{proof theory!ordinal
  analysis} of a theory~$T$, one typically produces an ordinal notation system
for ordinals less than some ordinal $\alpha$ such that for every $\beta <
\alpha$, the formalization~$\mathit{TI}(\beta)$ of the transfinite
induction\index{induction!transfinite} principle for~$\beta$ is provable
in~$T$. In practice, using transfinite induction up to~$\alpha$ itself and
otherwise only strictly finitary methods, one can prove the consistency
of~$T$.  The fact that induction up to $\epsilon_0$ establishes the
consistency of $\mathit{PA}$, together with the result of \cite{Gentzen:43}
that shows that for all $\beta < \epsilon_0$, $\mathit{PA}$ proves
$\mathit{TI}(\beta)$ for all $\beta < \epsilon_0$ constitutes an ordinal
analysis of~$\mathit{PA}$, and we say that $\epsilon_0$ is the proof theoretic
ordinal of~$\mathit{PA}$.

Proof theory in the tradition of Gentzen\index{Gentzen, G.} and
Sch\"utte\index{Sch\"utte, K.} as well as Takeuti\index{Takeuti, G.} has
focused on such ordinal analyses of theories of increasing strength.  In
recent work, Rathjen\index{Rathjen, M.}  \shortcite{Rathjen:05a,Rathjen:05b}
has pushed the boundaries of this approach in giving an ordinal analysis of a
very strong subsystem of analysis called $\Pi_2^1$-comprehension.\footnote{See
  \cite{Pohlers:87}\index{Pohlers, W.} for a survey of the work in the
  Sch\"utte school, and \cite{Pohlers:98} for a more recent technical survey.}
The consistency proofs in this tradition are, for the most part, based on the
approach of \cite{Schutte:60}\index{Sch\"utte, K.}, which uses a variant of
Gentzen\index{Gentzen, G.}'s formalization using infinitary derivations.  A
second tradition has pursued ordinal analysis using extensions of
Ackermann\index{Ackermann, W.}'s $\epsilon$-substitution method\index{epsilon
  substitution@$\epsilon$-substitution method}
\shortcite{Ackermann:40}\index{Ackermann, W.}, for examples see
\cite{Mints:99}\index{Mints, G.}\index{Tupailo, S.}  and
\cite{Arai:03}\index{Arai, T.}.

Although generalized Hilbert programs in this tradition have certainly
produced important mathematical work, its philosophical underpinnings are
thin. \cite{Takeuti:87} attempted to give a finitary justification for the
proof theoretic ordinal~$\epsilon_0$\index{epsilon zero@$\epsilon_0$}, but it
is unclear to what extent more complex ordinal notation systems are finitarily
acceptable.  Even if one concedes, as, e.g., Sch\"utte\index{Sch\"utte, K.}
does, that the consistency proofs\index{consistency proof} in question are
constructive (although no longer strictly finitary), it is still unclear what
the philosophical significance of the technical results is.
\cite[366]{Feferman:88a}\index{Feferman, S.}  offers this assessment:
\begin{quote}
  [A]s the systems of ordinal notation used for consistency proofs of stronger
  and stronger theories become more and more complicated, the significance to
  noncognoscenti of what is thereby accomplished decreases in inverse
  proportion. Thus, on the one hand, to say that one has obtained a
  constructive consistency proof\index{consistency proof} of a
  theory~$T$---without saying anything more---is too general to be
  informative; and, on the other hand, to say that the proof has been carried
  out by transfinite induction on a certain complicated recursive ordering for
  some very large ordinal tells us nothing about what constructive principles
  are involved in the proof of this well-ordering.\footnote{For a more
    forceful criticism of proof theory in this tradition, see
    \cite{Kreisel:76}\index{Kreisel, G.}.}
\end{quote}

Another important proof-theoretical approach in which the analysis of systems
of classical mathematics is accomplished using a generalization of the
finitary standpoint is that of functional interpretations\index{functional
  interpretation|(}.  The model for this approach is G\"odel\index{G\"odel,
  K.}'s \emph{Dialectica} interpretation\index{Dialectica interpretation|(}
\shortcite{Godel:58}.  The \emph{Dialectica} interpretation shows how one can
reduce an infinitary theory~$T$ (in this case, intuitionistic first-order
arithmetic) to a quantifier-free theory~$F$.\footnote{Via the interpretation
  of classical arithmetic in intuitionistic arithmetic
  \cite{Gentzen:33,Godel:33b}\index{Gentzen, G.}, the \emph{Dialectica}
  interpretation also yields a functional interpretation of classical
  arithmetic.}  An ordinal analysis of a theory does something similar, for
instance, one can view Gentzen\index{Gentzen, G.}'s consistency
proof\index{consistency proof} as reducing Peano arithmetic\index{Peano
  arithmetic} to a quantifier-free theory (\textit{PRA}\index{primitive
  recursive arithmetic}) extended by a certain infinitary induction principle
($\mathit{TI}(\epsilon_0)$).  In the case of functional interpretations, the
quantifier-free theory~$F$ is also not strictly finitary: it does not
just mention finite objects but also certain infinitary objects, viz.,
functionals of finite type.  A functional interpretation can be seen as a
reduction of the infinitary theory~$T$ to the theory of functionals~$F$ in
question.  The approach using functional interpretations has the following
advantage over the Gentzen\index{Gentzen, G.}-Sch\"utte approach.  It is a
consequence of the reduction of $T$ to $F$ that every recursive function which
can be proved to be total in~$T$ is represented by a term of~$F$.  Because the
functionals of~$F$ in practice characterize natural classes of functions, a
functional interpretation yields an appealing analysis of the computational
content of~$F$. Moreover, the conceptual import of the reduction is more
apparent than in the case of ordinal analysis: already in the case of
$\mathit{PA}$, G\"odel\index{G\"odel, K.}'s functionals of finite type provide
a much clearer account of the character of the constructive methods appealed
to than induction up to~$\epsilon_0$.\footnote{For an excellent survey of this
  approach, see \cite{Avigad:98}\index{Avigad, J.}\index{Feferman,
    S.}.}\index{Dialectica interpretation|)}\index{functional
  interpretation|)}

\subsection{Relativized Hilbert programs}\label{relativized-hp}

A philosophically more satisfactory continuation of Hilbert's program in proof
theoretic terms has been suggested by Kreisel\index{Kreisel, G.}
\shortcite{Kreisel:54,Kreisel:68,Kreisel:83} and has been elaborated
especially by Feferman\index{Feferman, S.}.  This work proceeds from a wider
conception of Hilbert's program as an attempt to justify ideal mathematics by
restricted means.  On this conception, the aim of Hilbert's proof theory was
to show that, at least as far as a certain class of real propositions is
concerned, ideal mathematics does not go beyond real mathematics, and in this
sense finitary mathematics is a foundation for ideal mathematics. A finitary
consistency proof\index{consistency proof} of the kind envisaged by
Hilbert\index{Hilbert, D.} would have accomplished this for all of classical
mathematics.

The scope of the foundational project, however, need not necessarily be all of
higher mathematics.  So-called relativized Hilbert programs are projects in
which one considers certain fragments of higher mathematics as the theory for
which a foundation is sought (and indeed, also theories stronger than finitism
as candidates for the reducing theory, e.g., predicative theories).  Examples
of these are Feferman\index{Feferman, S.}'s work on explicit
mathematics\index{explicit mathematics} and predicative subsystems of
analysis, and to some extent also the Friedman-Simpson program of reverse
mathematics (see below).  What is common to these approaches to mathematical
foundations is that they concentrate on so-called proof-theoretic reductions
of systems of classical mathematics to more restricted systems.  The reduction
is carried out using finitist means, and therein lies its philosophical
significance.
 
A foundational reduction\index{foundational reduction|(} in
Feferman's\index{Feferman, S.|(} sense \shortcite{Feferman:88a,Feferman:93} is
accomplished if it can be shown that a body of mathematics which is justified
by a foundational framework ${\cal F}_1$ (e.g, finitary, constructive,
predicative, infinitary, set-theoretic) can already be justified, in a certain
sense, in a weaker, or stricter foundational framework ${\cal F}_2$.  This is
in general not possible in a wholesale fashion, however, partial foundational
reductions can and have been achieved.  Suppose a theory $T_1$ is justified by
a foundational framework ${\cal F}_1$, and a theory $T_2$ by a weaker
framework ${\cal F}_2$. A proof theoretic reduction\index{proof theory!proof
  theoretic reduction} of $T_1$ to $T_2$ (conservative for $\Phi$) is a
partial recursive function $f$ such that
\begin{enumerate}
\item Whenever $x$ is (the code of) a proof in $T_1$ of a formula
(with code) $y$ in $\Phi$, then $f(x)$ is (the code of) a proof of $y$
in $T_2$, and
\item $T_2$ proves the formalization of (1).
\end{enumerate}
If there is such a function $f$, we write $T_1 \leq T_2 [\Phi]$.  Now if $T_1$
is directly justified by a foundational framework ${\cal F}_1$, and $T_2$ by
${\cal F}_2$, then, so Feferman\index{Feferman, S.}, a proof-theoretic
reduction that establishes $T_1 \leq T_2 [\Phi]$ is a partial foundational
reduction of ${\cal F}_1$ to ${\cal F}_2$.  Clause (2) in the definition
ensures that the reduction (the function~$f$) itself is justified by the
weaker framework ${\cal F}_2$.  In the reductions achieved in practice, it
turns out that $f$ is actually primitive recursive and the formalization of
(1) can even be proved in primitive recursive arithmetic
$\mathit{PRA}$\index{primitive recursive arithmetic}. Since
$\mathit{PRA}$\index{primitive recursive arithmetic} is directly justified by
the finitary framework, such partial foundational reductions are therefore all
finitarily justified.  Feferman's\index{Feferman, S.} main philosophical
conclusion from the possibility of giving such foundational reductions is
this: The main argument for set-theoretical realism is the
Quine-Putnam\index{Quine, W. V.  O.}\index{Putnam, H.} indispensability
argument\index{indispensability argument}, which proceeds from the premises
that set-theory is indispensable to science.  Feferman has shown, first, that
much, if not all, of scientifically applicable mathematics can actually be
formalized in much weaker systems (e.g., Feferman's system $W$, which is
justified by a predicative foundational framework), and second, that
predicative mathematics can be reduced to the countably infinite (in the sense
that there is a partial foundational reduction of predicative mathematics to
countably infinite mathematics, given by a proof-theoretic reduction of $W$ to
Peano Arithmetic~\textit{PA})\index{Peano arithmetic}.  He concludes that,
\begin{quote}
  even if one accepts the indispensability argument, practically nothing
  philosophically definite can be said of the entities which are then supposed
  to have the same status---ontologically and epistemologically---as the
  entities of natural science.  That being the case, what do the
  indispensability arguments amount to? As far as I'm concerned, they are
  completely vitiated. \cite{Feferman:93a}
\end{quote} 
Independently of the question of mathematical realism and of the scope and
force of the indispensability arguments, proof-theoretic reductions give
precise answers to questions of the relation between foundational frameworks.
Since a proof-theoretic reduction of $T_1$ to $T_2$ also yields a consistency
proof\index{consistency proof} of $T_1$ in $T_2$ (i.e., a relative consistency
result), establishing a proof-theoretic reduction also provides a solution to
Hilbert's program relativized to $T_1$ and $T_2$.  Feferman summarizes the
importance of proof-theoretic reductions thus:
\begin{quote}
  In general, the kinds of results presented here serve to sharpen what is to
  be said in favor of, or in opposition to, the various philosophies of
  mathematics such as finitism, predicativism, constructivism, and
  set-theoretical realism. Whether or not one takes one or another of these
  philosophies seriously for ontological and/or epistemological reasons, it is
  important to know which parts of mathematics are in the end justifiable on
  the basis of the respective philosophies and which are not.  The uninformed
  common view---that adopting one of the non-platonistic positions means
  pretty much giving up mathematics as we know it---needs to be drastically
  corrected, and that should also no longer serve as the last-ditch stand of
  set-theoretical realism.  On the other hand, would-be nonplatonists must
  recognize the now clearly marked sacrifices required by such a commitment
  and should have well-thought out reasons for making them. \cite{Feferman:93}
\end{quote}
Proof theorists have obtained a number of such results, including reductions
of theories which on their face require a significant amount of ideal
mathematics for their justification (e.g., subsystems of analysis) to finitary
systems.\footnote{For a discussion of the philosophical significance of such
  proof theoretic reductions, see \cite{Feferman:00} and
  \cite{Hofweber:00}\index{Hofweber, T.}.}\index{Feferman,
  S.|)}\index{foundational reduction|)}

The program of so-called reverse mathematics\index{proof theory!reverse
  mathematics}\index{reverse mathemtaics}
developed by, in particular, Friedman\index{Friedman, H.} and
Simpson\index{Simpson, S.}, is another continuation of Hilbert's program.  In
the face of G\"odel\index{G\"odel, K.}'s results showing that not all of
classical mathematics can be reduced to the finitary, they seek to answer the
question: how much of classical mathematics can be so reduced?  Reverse
mathematics aims to give a precise answer to this question by investigating
which theorems of classical mathematics are provable in weak subsystems of
analysis which are reducible to finitary mathematics (in the sense discussed
above).  A typical result is that the Hahn-Banach theorem
of functional analysis is provable in a theory known as $WKL_0$ (for ``weak
K\"onig lemma''); $WKL_0$ is proof-theoretically reducible to
$PRA$\index{primitive recursive arithmetic} for $\Pi^0_2$ sentences (i.e.,
sentences of the form $\forall x\exists y\, A(x, y)$.\footnote{See
  \cite{Simpson:88} for an overview, \cite{Simpson:99} for a technical
  introduction to reverse mathematics, and also the collection
  \cite{Simpson:05}.}

Reverse mathematics in this tradition is primarily concerned with infinitary
theories, e.g., subsystems of analysis.  G\"odel's theorems show, however,
that not even all truths of first-order number theory are provable in Peano
arithmetic\index{Peano arithmetic}, and hence that not even the domain of all
arithmetical truths can be given a foundation on finitary principles.  This
suggests the questions of whether there are ``mathematically interesting''
statements of number theory which are not provable in systems that can be
reduced to the finitary. The most celebrated result in this regard is the
proof by \cite{Paris:77}\index{Paris, J.}\index{Harrington, L.} that a version
of the finite Ramsey theorem is not provable in Peano arithmetic\index{Peano
  arithmetic}.  However, this and other examples of independent number
theoretic statements are constructed specifically to be independent of Peano
arithmetic.  It turns out that a great many ``ordinary'' number theoretic
results are provable even in weak fragments of first-order number theory, and
this has led Friedman to conjecture that ``every theorem published in the
\emph{Annals of Mathematics} whose statement involves only finitary
mathematical objects (i.e., what logicians call an arithmetical statement) can
be proved in elementary arithmetic.''  (Here, elementary arithmetic is a very
weak theory which can be proved consistent by primitive recursive methods.)
\cite{Avigad:03}\index{Avigad, J.} gives an excellent survey of the issues and
results related to this conjecture and places it in the philosophical context
of Hilbert's program.\footnote{See also \cite{Raatikainen:03} on the current
  status of the various branches of proof-theoretic research relating to
  Hilbert's program.}

The results surveyed in theis section are all natural continuations of
Hilbert's\index{Hilbert, D.} original ideas.  A central aspect recent
proof-theoretical investigations and of Hilbert's original program alike that
they study formalized systems using meta-mathematical tools with the aim of
understanding the structure and content of these systems.  Hilbert's original
consistency project, the conservativity project that Kreisel\index{Kreisel,
  G.} and others interpret Hilbert as having been engaged in, as well as
reductive proof theory are all examples of this, and this is also part of the
reason why many practicing proof theorists see themselves as still working on
Hilbert's program. Ordinal analysis, functional interpretations, proof
theoretic reductions and reverse mathematics are only some of the most
prominent areas of proof theory, and those most explicitly situated in the
tradition of Hilbert's program.  Many other areas of proof theory other than
those directly concerned with consistency and foundational reductions of
theories are related to the aims of Hilbert's program, e.g., the
no-counterexample interpretation\index{no-counterexample interpretation}
\cite{Kreisel:51,Tait:05a}\index{Tait, W. W.} and work on the structure and
complexity of formal proofs\index{proof theory!complexity of proofs} more
generally \cite{Pudlak:98}\index{Pudl\'ak, P.}.

\section{Conclusion}

Although it has been a commonplace in the literature on the philosophy of
mathematics and logic from 1950 onward that Hilbert's program has not only
been ``killed'' by G\"odel\index{G\"odel, K.}'s incompleteness
theorems\index{incompleteness theorems} but that it was over-ambitions if not
ill-conceived from the start, in the current literature a more positive
evaluation has emerged.  This is in large part due to the attention which
unpublished writings in the Hilbert school (especially lecture notes to
Hilbert's\index{Hilbert, D.} courses) have received recently, as well as to
the availability of more of the published writings in English translation
(e.g., in \cite{Ewald:96} and \cite{Mancosu:98}).  But it is also due to a
growing recognition that the common characterizations of Hilbert's program are
a caricature, and to a clearer philosophical engagement with the recent
results of proof theory.  For instance, Ramsey's characterization that,
according to Hilbert, ``Mathematics proper is thus regarded as a sort of game,
played with meaningless marks on paper rather like noughts and crosses''
\cite[231]{Ramsey:26}\index{Ramsey, F. P.}, and the view that Hilbert held a
naive formalist and instrumentalist view of mathematics have been criticized
by various writers.  It remains to be seen in what ways Hilbert's
philosophical views can be resuscitated (in the manner in which, e.g.,
Frege\index{Frege, G.}'s logicist program has experienced a renaissance in the
neo-logicist writings of, e.g., Boolos\index{Boolos, G.}, Heck\index{Heck,
  R.}, and Hale\index{Hale, B.} and Wright\index{Wright, C.}).  It should be
clear in any case from the discussion in the preceding section that ideas
closely related to Hilbert's own \emph{have} been hugely successful.  And it
is also clear from the recent historical studies of Hilbert's unpublished
writings as well as from the study of the proof theoretical practice of the
Hilbert school that the ideas which underpin much of recent and current
proof theoretical research are not merely ``inspired by'' Hilbert's program.
Hilbert's fundamental aim was, all along, to make mathematical reasoning
amenable to mathematical investigation, and to carry out such an investigation
which yields an analysis of non-constructive reasoning in terms of restricted
methods.  Hilbert, of course, emphasized consistency of non-constructive
systems as the most interesting property to be investigated, and emphasized
finitary methods as those in terms of which such an analysis should be carried
out.  But even in the 1920s, in the practice of consistency
proofs\index{consistency proof} in the work of Ackermann\index{Ackermann, W.},
Bernays\index{Bernays, P.}, and von Neumann\index{von Neumann, J.}, among
others, more broadly constructive methods were employed in this analysis, and
results about properties other than consistency were obtained.
Gentzen\index{Gentzen, G.}'s work of the 1930s and subsequent proof
theoretical studies should thus not be seen as merely a response to
G\"odel\index{G\"odel, K.}'s incompleteness results\index{incompleteness
  theorems}, but more broadly as advancing Hilbert's original aim of
investigating the structure of mathematical reasoning.  Seen in this light,
again, proof theory as a foundational enterprise is very much alive.  Although
G\"odel\index{G\"odel, K.}'s theorems show that Hilbert's original
expectations about what exactly can be analyzed in which way and with what
restricted methods can not be fulfilled, proof theory and Hilbert's aims more
generally have been advanced tremendously over the last
half-century.\index{Hilbert's program|)}

\subsection*{Acknowledgements}

This paper is a substantially revised and expanded version of my entry in the
\emph{Stanford Encyclopedia of Philosophy} \cite{Zach:03a}.  I would like to
acknowledge the support of the Social Sciences and Humanities Research Council
of Canada, as well as helpful comments from Jeremy Avigad, Aleksandar
Ignjatovi\v{c}, Paolo Mancosu, and William Tait.

\end{document}